\documentclass[journal]{IEEEtran}
\usepackage{amsmath,amsfonts}
\usepackage{amsthm}
\usepackage{algorithm,algorithmic}
\usepackage{array}
\usepackage[caption=false,font=normalsize,labelfont=sf,textfont=sf]{subfig}
\usepackage{textcomp}
\usepackage{stfloats}
\usepackage{url}
\usepackage{verbatim}
\usepackage{graphicx}
\usepackage{cite}
\usepackage{enumitem}
\def\BibTeX{{\textup B\kern-.05em{\sc i\kern-.025em b}\kern-.08em
    T\kern-.1667em\lower.7ex\hbox{E}\kern-.125emX}}
\usepackage{balance}

\usepackage{amssymb}

\newtheoremstyle{indented}
  {0pt} 
  {0pt} 
  {\itshape} 
  {\parindent} 
  {\itshape} 
  {:}
  { } 
  {}

\theoremstyle{indented}
\newtheorem{theorem}{Theorem}
\newtheorem{proposition}{Proposition}
\newtheorem{definition}{Definition}

\newtheorem{fact}{Fact}

\DeclareMathOperator{\sign}{sign}

\begin{document}
\title{
Sparse Regularization by Smooth Non-separable Non-convex Penalty Function Based on Ultra-discretization Formula
}
\author{Natsuki Akaishi, \IEEEmembership{Student Member, IEEE}, Koki Yamada, \IEEEmembership{Member, IEEE}, Kohei Yatabe, \IEEEmembership{Member, IEEE}
\thanks{Manuscript received XXXX XX, XXXX; revised XXXXX XX, XXXX; accepted XXXXX XX, XXXX. Date of publication XXXXX XX, XXXX; date of current version XXXXX XX, XXXX. The associate editor was XXXXX XXXX. \\
The authors are with Tokyo University of Agriculture and Technology,
Tokyo 184-8588, Japan (e-mail: natsu61aka@gmail.com; k-yamada@go.tuat.ac.jp; yatabe@go.tuat.ac.jp).}}

\markboth{IEEE TRANSACTIONS ON SIGNAL PROCESSING,~Vol.~XX, No.~X, XXXX~20XX}%
{AKAISHI \MakeLowercase{\textit{et. al.}}: Sparse Regularization by Smooth Non-separable Non-convex Penalty Function Based on Ultra-discretization Formula}

\maketitle

\begin{abstract}

In sparse optimization, the $\ell_{1}$ norm is widely adopted for its convexity, yet it often yields solutions with smaller magnitudes than expected.
To mitigate this drawback, various non-convex sparse penalties have been proposed.
Some employ non-separability, with ordered weighting as an effective example, to retain large components while suppressing small ones.
Motivated by these approaches, we propose ULPENS, a non-convex, non-separable sparsity-inducing penalty function that enables control over the suppression of elements.
Derived from the ultra-discretization formula, ULPENS can continuously interpolate between the $\ell_{1}$ norm and a non-convex selective suppressing function by adjusting parameters inherent to the formula.
With the formula, ULPENS is smooth, allowing the use of efficient gradient-based optimization algorithms.
We establish key theoretical properties of ULPENS and demonstrate its practical effectiveness through numerical experiments.

\end{abstract}

\begin{IEEEkeywords}
Sparse signal estimation, non-convex regularization, Lipschitz constant, smooth approximation, quasi-Newton method.
\end{IEEEkeywords}

\section{Introduction}
\label{sec:intro}

Sparse representations play a central role in various applications, including noise reduction, image deblurring, and compressed sensing \cite{zhang2021overview,gao2020spectral,wright2010sparse,milton2020controller,TIAN2022146,chen2021solving,Marques2019}.
To address the ill-conditioned or severely underdetermined system of linear equations arising in those applications, algorithms commonly aim to identify solutions that are sparse or approximately sparse.
A widely used approach to find a sparse approximate solution $\mathbf{x}\in\mathbb{R}^{N}$ is to minimize the following cost function $f: \mathbb{R}^{N} \to \mathbb{R}$,
\begin{equation}
\label{eq:cost}
    f(\mathbf{x}) = g_{\mathbf{s}}(\mathbf{x}) + \gamma\,r(\mathbf{x}),
\end{equation}
where $\mathbf{s}\in\mathbb{R}^{M}$ is an observation, $g_{\mathbf{s}}(\mathbf{x})$ is a data fidelity term, $\gamma\,r(\mathbf{x})$ is a sparse regularization (or penalty) term, and $\gamma > 0$ is a parameter that balances the two terms.
A typical example is $g_{\mathbf{s}}(\mathbf{x}) = \frac{1}{2}\|\mathbf{A}\mathbf{x} - \mathbf{s}\|^2_{2}$, where $\mathbf{A}\in\mathbb{R}^{M\times N}$ is an observation matrix, and $r(\mathbf{x}) = \|\mathbf{x}\|_1$, i.e., the $\ell_{1}$ norm.
A wide range of applications use the $\ell_{1}$ norm for sparse regularization \cite{l1norm}.

\begin{figure}
    \centering
    \includegraphics[width=0.85\linewidth]{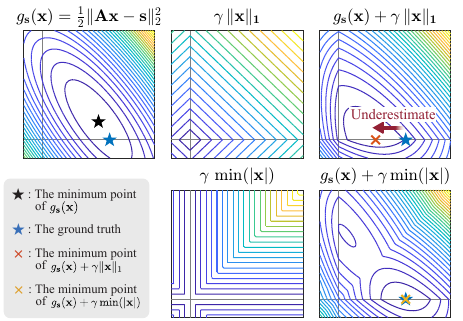}
    \caption{Contour plots of $g_{\mathbf{s}}=\frac{1}{2}\|\mathbf{A}\mathbf{x} - \mathbf{s}\|_{2}^2$ (top left), the ($\gamma$-scaled) $\ell_{1}$ norm (top center), the ($\gamma$-scaled) minimum absolute value function (bottom center), $\frac{1}{2}\|\mathbf{A}\mathbf{x} - \mathbf{s}\|_{2}^2 + \gamma\|\mathbf{x}\|_{1}$ (top tight) and $\frac{1}{2}\|\mathbf{A}\mathbf{x} - \mathbf{s}\|_{2}^2 + \gamma\min(|\mathbf{x}|)$ (bottom right).
    The colors of the contours indicate the function values: blue represents smaller values, while yellow represents larger values.
    The marks on the plot are as shown in the legend at the bottom left.
    The axes of each variable are shown in gray lines.}
    \label{fig:l1vsMin}
\end{figure}

While $\ell_{1}$-norm-based optimization is effective in promoting sparsity, it tends to underestimate high-amplitude components, which are often the most relevant parts of the signal.
An example of underestimation in the two-dimensional case is shown in the top boxes of Fig.~\ref{fig:l1vsMin}.
This example addresses the problem of estimating the ground truth (blue star), which is a sparse, nonzero solution that minimizes $g_{\mathbf{s}}$.
As shown in the upper right, the minimum point of the cost function (top right) lies on the axis.
This means that the solution (red cross) is sparse, but the amplitude of the solution is smaller than the ground truth.
This occurs because the value of the $\ell_{1}$ norm increases with the amplitude of each element, as shown in the top center.

Numerous studies have demonstrated that replacing the $\ell_{1}$-norm with a properly selected non-convex function can lead to improved performance \cite{Zhang2010nearly,yin2015minimization,LYU2013413,candes2008enhancing,lanza2019sparsity,Ivan2017,Bogdan2015,Tipping2001,Wipf2011,tao2022minimization,wipf2010iterative,Ivan2017nonsep,zeng2014decreasing,sander2023fast,sasaki2024}
In particular, the importance of non-separable penalties, which account for the relationships among elements, has been recognized \cite{sander2023fast,wipf2010iterative,Ivan2017nonsep,zeng2014decreasing,sasaki2024}.
Some of these penalties incorporate the magnitude relationships among elements by sorting elements and applying selective weights \cite{zeng2014decreasing,sander2023fast,sasaki2024}.
In \cite{sasaki2024}, decreasing weights are applied, with no weight assigned to the larger elements.
This allows small elements to be suppressed to zero while preserving the larger ones.
As the extreme example, consider the two-dimensional case, setting the smaller weight in \cite{sasaki2024} to zero corresponds to the minimum absolute value function.
With the minimum absolute value function, as shown in the bottom right of Fig.~\ref{fig:l1vsMin}, the ground truth and the minimum point of the cost function (yellow cross) are nearly identical.
However, its straightforward extension to higher dimensions is non-trivial.

In this paper, we propose ULPENS, a non-convex, non-separable and smooth penalty function inducing sparsity, designed to overcome the underestimation issue of the $\ell_{1}$ norm.
ULPENS is constructed from an approximation of the minimum absolute value function, derived from the ultra-discretization formula. 
It provides a smooth interpolation between the $\ell_{1}$ norm and a selective suppression function that attenuates smaller components.
Moreover, as a by-product, ULPENS is smooth.
Its smoothness eliminates the need for non-smooth optimization algorithms, which are a common requirement in many existing sparse penalties, and enables the use of faster optimization methods such as the quasi-Newton method.
We establish key theoretical properties of ULPENS, including bounds on the Lipschitz constant of its gradient, and demonstrate its practical effectiveness through numerical experiments.

The rest of the paper is organized as follows.
Section II reviews optimization algorithms and conventional sparsity-inducing functions. 
Section III presents the proposed penalty, ULPENS, along with the ultra-discretization formula.
Section IV discusses its properties and implementations. 
Section V provides the experiments validating the effectiveness of the proposed method. 
Section VI concludes this paper.

\section{Preliminaries}

\subsection{Optimization Algorithms}
\label{ssec:smtOpt}

Optimization problems, such as minimization of \eqref{eq:cost}, are solved by numerical algorithms.
As is common in sparse optimization, we assume that $g_{\mathbf{s}}$ in \eqref{eq:cost} is smooth, i.e., it has a continuous gradient.  
In such cases, the choice of algorithms for minimizing \eqref{eq:cost} depends on the smoothness of the regularization function $r$, which is usually non-smooth in sparse optimization.
To address non-smooth problems, the proximal gradient method (PGM) \cite{bauschke2017}, the alternating direction method of multipliers (ADMM) \cite{boyd2004convex}, and the primal-dual splitting method (PDS) \cite{playing} can be used.
If $r$ is smooth, gradient-based approaches, such as the gradient descent method (GD), can be applied.
Various techniques to accelerate these algorithms are also widely used.
Common acceleration techniques include the heavy ball method \cite{POLYAK19641,le2024nonsmooth,Aujol2022}, Nesterov acceleration \cite{nesterov2018lectures, nesterov2005smooth,blumensath2008iterative,buccini2020general}, and the quasi-Newton method \cite{broyden1967quasi,pmlr-v162-jin22b,berahas2022quasi}, where the quasi-Newton method is only applicable to smooth optimization.

\subsection{$\ell_{1}$ Norm Minimization and Underestimation Problem}
\label{ssec:SIP}

In sparse optimization, instead of the $\ell_{0}$ norm%
\footnote{
Although the ``$\ell_{0}$ norm'' is not a norm, we use the notation $\|\cdot\|_{0}$.}, the $\ell_{1}$ norm and its extensions are often used.
The $\ell_{1}$ norm is defined as
\begin{equation}
    \|\mathbf{x}\|_{1} = \sum^{N}_{n=1}|x_n|,
\end{equation}
where $\mathbf{x}\in\mathbb{R}^{N}$.
Note that the derivative of the absolute value function is ill-defined at the origin.
To facilitate the later discussion on the relationship between the proposed method and the $\ell_{1}$ norm, we treat the derivative of the absolute value function at the origin as 0, making it equivalent to the signum function.
This also makes the gradient of the $\ell_{1}$ norm equivalent to the element-wise signum function.

The $\ell_{2,1}$ norm is one of the extensions of the $\ell_{1}$ norm to induce a structured sparsity \cite{L21}.
The $\ell_{2,1}$-norm is defined as
\begin{equation}
\label{eq:L21}
    \|\mathbf{x}\|^{\mathcal{G}}_{2,1} = \|\left[\|\mathbf{x}_{\mathcal{G}_1}\|_{2}, \dots, \|\mathbf{x}_{\mathcal{G}_M}\|_{2}\right]^{\mathsf{T}}\|_{1},
\end{equation}
where $\|\cdot\|_{2}$ is the $\ell_{2}$ norm, $(\cdot)^{\mathsf{T}}$ denotes transpose, $\mathcal{G}_{1}, \dots, \mathcal{G}_{M}$ ($1 < M < N$) are index sets that satisfy $\mathcal{G}_{k} \cap \mathcal{G}_{l} = \emptyset$ for any $k \neq l$, and $\bigcup_{m=1}^{M} \mathcal{G}_{m} = \{1, \dots, N\}$.
$\mathbf{x}_{\mathcal{G}_{m}} \in \mathbb{R}^{I_m}$ is a vector in $\mathcal{G}_{m}$, where $I_m$ is the number of elements of $\mathcal{G}_{m}$.

Although the $\ell_{1}$ norm minimization is widely used, it often results in a solution whose magnitude is smaller than the ground truth because the $\ell_{1}$ norm penalizes all entries of the variable.
The same issue occurs with the $\ell_{2,1}$ norm.
In this paper, we call this problem as the underestimation problem.

\subsection{Non-convex Non-separable Sparse Penalties}
\label{ssec:nonsep}

\begin{figure}
    \centering
    \includegraphics[width=0.78\linewidth]{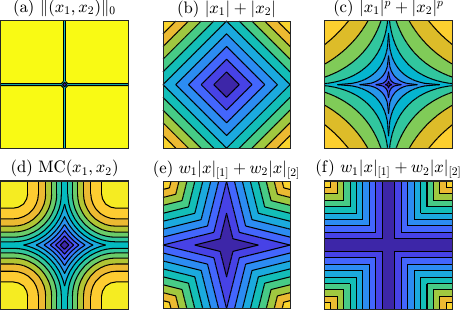}
    \caption{Examples of contour plots for conventional sparse regularization functions in the two-dimensional space.
    The origins are at the center, with the vertical and horizontal axes representing each variable.
    The color indicates the function values: blue represents smaller values, while yellow represents larger values.
    (a) is the $\ell_{0}$ norm, (b) is the $\ell_{1}$ norm, (c) is the $p$-powered $\ell_{p}$ norm with $p = 0.5$, (d) is the MC penalty, (e) is the OWL penalty with $w_1=0.8, w_2=0.2$, and (f) is the OWL penalty with $w_1=0$, $w_2=1$.
    }
    \label{fig:convSparse}
\end{figure}

To address the underestimation problem in sparse optimization, many non-convex penalties have been proposed \cite{Ivan2017,wipf2010iterative,Zhang2010nearly,yin2015minimization,zeng2014decreasing,tao2022minimization,candes2008enhancing,LYU2013413,sasaki2024,Bogdan2015,Tipping2001,Wipf2011,Ivan2017nonsep,lanza2019sparsity,sander2023fast}.
These penalties can be categorized into separable penalties \cite{Zhang2010nearly,yin2015minimization,LYU2013413,candes2008enhancing} and non-separable penalties \cite{tao2022minimization,wipf2010iterative,Ivan2017,Tipping2001,Wipf2011,Ivan2017nonsep,lanza2019sparsity,Bogdan2015,zeng2014decreasing,sander2023fast,sasaki2024}.
In particular, the importance of non-separability has been recognized \cite{Wipf2011,Ivan2017nonsep}.
This is because non-separable penalties allow a selective regularization of input elements.

Here, we introduce three non-convex penalties \cite{Ivan2017,LYU2013413,sasaki2024,zeng2014decreasing} shown in Figs.~\ref{fig:convSparse}(c) to (f), with a particular focus on the non-separable penalties in (e) and (f).
The $p$-powered $\ell_{p}$ penalty \cite{LYU2013413} and the minimax concave (MC) penalty \cite{Zhang2010nearly} in Figs.~\ref{fig:convSparse}(c) and (d), respectively, are separable.
A non-separable extension of the MC penalty, the generalized MC (GMC) penalty, has also been proposed \cite{Ivan2017}.
Another example of a non-separable penalty is the ordered weighted $\ell_{1}$ (OWL) penalty \cite{zeng2014decreasing,sasaki2024}, which is defined as
\begin{equation}
\label{eq:WL1}
    s_{\mathbf{w}}(\mathbf{x}) = \sum^{N}_{n=1} w_n|x|_{[n]},
\end{equation}
where $|x|_{[n]}$ is the $n$th largest value of $|x_1|,\dots,|x_N|$ i.e., $|x|_{[1]}\geq\dots\geq|x|_{[N]}$, and $w_n\in[0,+\infty)$ is a weight.
Fig.~\ref{fig:convSparse}(e) and (f) show examples of the contour plot of OWL, where the weights are assigned such that $ w_1 = 0.2$ and $w_2 =0.8$ in (e) and are $w_1 = 0$, $w_2 = 1$ in (f).
Note that the latter is equivalent to the minimum absolute value function, $\min(|x_1|,|x_2|)$.
This function can selectively suppress elements with small amplitudes, resulting in the avoidance of the underestimation problem, as illustrated in the bottom boxes in Fig.~\ref{fig:l1vsMin}.

\section{Proposed Method}
\label{sec:proposal}

\begin{figure}
    \centering
    \includegraphics[width=1\linewidth]{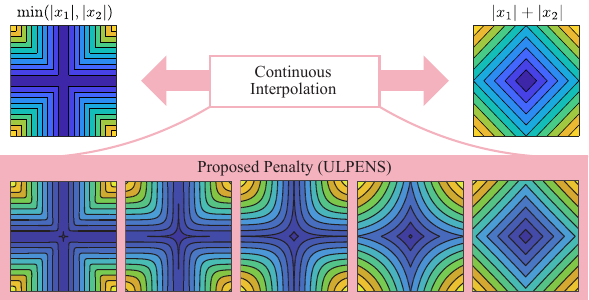}
    \caption{
    Concept of the proposed penalty (ULPENS) and its contour plots in the two-dimensional space.
    The origins are at the center.
    The top left box is the minimum absolute value function, the top right box is the $\ell_{1}$ norm, and the bottom boxes are ULPENS.
    }
    \label{fig:pena}
\end{figure}

In this section, we propose a smooth, non-convex sparsity-inducing penalty, named ULPENS (\underline{UL}tra-discretization-formula-based sparsity-inducing \underline{PE}nalty with \underline{N}on-convexity and \underline{S}moothness), which is designed to prevent the underestimation of non-zero elements by selectively suppressing smaller elements.
Let us first illustrate the contour plot of ULPENS in the two-dimensional space in Fig.~\ref{fig:pena}.
As shown, ULPENS continuously interpolates between the $\ell_{1}$ norm and the following minimum absolute value function: 
\begin{equation}
\label{eq:propId}
    \min(|\mathbf{x}|)\quad(= \min(|x_1|,\dots,|x_N|)),
\end{equation}
where $|\cdot|$ denotes the element-wise absolute value.
The interpolation is achieved through an approximation inspired by the ultra-discretization formula.
This allows continuous adjustment between the $\ell_1$ norm and the minimum absolute value function, yielding a penalty function that can effectively control the amount of suppression of larger elements.
Moreover, as a by-product of the approximation, ULPENS allows the use of a fast smooth optimization solver.
We also propose an enhancement of the $\ell_{2,1}$ norm using ULPENS.

\subsection{Ultra-discretization Formula}

The following fact provides the ultra-discretization formula \cite{Nakazono2017} in \eqref{eq:max}, which is the key to the proposed method.

\begin{fact}
\label{fact:max}
Let $\mathbf{x}\in\mathbb{R}^{N}$, and $\lambda > 0$.
When $\lambda \to +0$, the following holds:
\begin{equation}
\label{eq:max}
    \lambda\ln\left(\sum^{N}_{n=1}\textup{e}^{\frac{x_n}{\lambda}}\right) \to \max\left(\mathbf{x}\right).
\end{equation}
\end{fact}
\indent\textit{Proof:} 
Let $x_{\max} = \max\left(\mathbf{x}\right)$.
The left-hand side of \eqref{eq:max} can be rearranged as
\begin{align*}
   \lambda\ln\left(\sum^{N}_{n=1}\textup{e}^{\frac{x_n}{\lambda}}\right)&= \lambda\ln\left(\textup{e}^{\frac{x_{\max}}{\lambda}}\sum^{N}_{n=1}\textup{e}^{\frac{x_n-x_{\max}}{\lambda}}\right) \\
     & = x_{\max} + \lambda\ln\left(\sum^{N}_{n=1}\textup{e}^{\frac{x_n-x_{\max}}{\lambda}}\right).
\end{align*}
The limit of the second term as $\lambda\to+0$ is given as follows. 
If $x_n \neq x_{\max}$, $\textup{e}^{{\frac{x_n-x_{\max}}{\lambda}}}\to +0$ because $x_n-x_{\max} < 0$ and $\lambda >0$.
If $x_n = x_{\max}$, $\textup{e}^{{\frac{x_n-x_{\max}}{\lambda}}}=1$.
Thus, $\sum^{N}_{n=1}\textup{e}^{\frac{x_n-x_{\max}}{\lambda}}\to +J$, where $J$ is the multiplicity of the maximum value.
Since $\ln\left(\sum^{N}_{n=1}\textup{e}^{\frac{x_n-x_{\max}}{\lambda}}\right)$ converges to the finite constant $\ln J$, $\lambda\ln\left(\sum^{N}_{n=1}\textup{e}^{\frac{x_n-x_{\max}}{\lambda}}\right)\to +0$ as $\lambda\to+0$.\qed

\subsection{Smoothly Approximated Absolute Value Function}
\label{ssec:smAbs}

The ultra-discretization formula in \eqref{eq:max} can smoothly approximate the absolute value function as follows.
\begin{definition}
\label{def:abs}
Let $\lambda > 0$.
A smoothly approximated absolute value function $h^{\lambda}_{\textup{abs}}: \mathbb{R} \to [\lambda\ln2,\infty)$ is defined as 
\begin{equation}
\label{eq:abs}
    h^{\lambda}_{\textup{abs}}(t) = \lambda\ln\left(\textup{e}^{\frac{t}{\lambda}} + \textup{e}^{-\frac{t}{\lambda}}\right). 
\end{equation}
\end{definition}
\begin{proposition}
\label{prop:abs}
$h^{\lambda}_{\textup{abs}}(t)$ converges to $|t|$ as $\lambda \to +0$, and is smooth.
Its derivative with respect to $t$ is given by the hyperbolic tangent function:
\begin{equation}
\label{eq:tanh}
    \frac{\textup{d}}{\textup{d}t} h^{\lambda}_{\textup{abs}}(t) = \frac{\textup{e}^{\frac{t}{\lambda}} - \textup{e}^{-\frac{t}{\lambda}}}{\textup{e}^{\frac{t}{\lambda}} + \textup{e}^{-\frac{t}{\lambda}}} = \tanh\left(\frac{t}{\lambda}\right).
\end{equation}

\end{proposition}
\indent\textit{Proof:} 
The convergence result follows from $|t| = \max(t,-t)$ and the application of \eqref{eq:max} to this expression.
The smoothness of $h^{\lambda}_{\textup{abs}}$ is ensured because it is a composite function of smooth functions, $\exp(\cdot)$ and $\ln(\cdot)$, with properly defined domains.
The derivative in \eqref{eq:tanh} is obtained by differentiating \eqref{eq:abs} with respect to $t$.
\qed

Hereafter, $\lambda$ is treated as the parameter associated with the approximation of the absolute value function.

\begin{figure}
    \centering
    \includegraphics[width=1\linewidth]{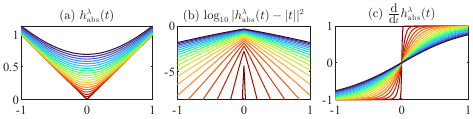}
    \caption{Examples of the smoothly approximated absolute value function.
    (a) shows $h^{\lambda}_{\textup{abs}}(t)$ in \eqref{eq:abs}.
    (b) shows the error $\log_{10}|h^{\lambda}_{\textup{abs}}(t) - |t||^2$.
    (c) shows $\frac{\text{d}}{\text{d}t}h^{\lambda}_{\textup{abs}}(t)$ in \eqref{eq:tanh}.
    The smaller the $\lambda$ in \eqref{eq:abs} and \eqref{eq:tanh}, the redder the color. 
    $\lambda$ changes linearly from 0.01 to 1.
    }
    \label{fig:approxAbs}
\end{figure}

The shape of the function $h^{\lambda}_{\textup{abs}}$ and its derivative are shown in Fig.~\ref{fig:approxAbs}.
The approximation eliminates the discontinuity in the derivative around the origin, as shown in Fig.~\ref{fig:approxAbs}(c).
As $\lambda$ decreases (redder), the approximation of the absolute value function becomes more accurate, as shown in Fig.~\ref{fig:approxAbs}(b).
The derivative $\frac{\textup{d}}{\textup{d}t}h^{\lambda}_{\textup{abs}}$ also converges to $\sign$, which is the derivative of the absolute value function, as shown in Fig.~\ref{fig:approxAbs}(c).

\subsection{Smoothly Approximated Minimum Value Function}
\label{ssec:smmin}

The ultra-discretization formula in \eqref{eq:max} can smoothly approximate the minimum value function as follows.
\begin{definition}
Let $\mu>0$.
A smoothly approximated minimum value function $h^{\mu}_{\text{min}}: \mathbb{R}^{N}\to\mathbb{R}$ is defined as
\begin{equation}
\label{eq:aprxMin}
    h^{\mu}_{\min}(\mathbf{x}) = -\mu\ln\left(\sum^{N}_{n=1}\textup{e}^{-\frac{x_n}{\mu}}\right).
\end{equation}
\end{definition}
\begin{proposition}
\label{prop:min}
$h^{\mu}_{\text{min}}(\mathbf{x})$ converges to $\min(\mathbf{x})$ as $\mu\to+0$, and is smooth.
Its gradient with respect to $\mathbf{x}$ is given by
\begin{equation}
\label{eq:aprxGradMin}
    [\nabla h^{\mu}_{\min}(\mathbf{x})]_{n} = \frac{\textup{e}^{-\frac{x_n}{\mu}}}{\sum^{N}_{n=1}\textup{e}^{-\frac{x_n}{\mu}}},
\end{equation}
where $[\cdot]_{n}$ denotes the $n$th element of the vector.
\end{proposition}
\indent\textit{Proof:}
Since $\min(\cdot) = -\max(-\cdot)$, $h^{\mu}_{\min}$ is derived by changing the sign of \eqref{eq:max}.
The convergence result is followed from \eqref{eq:max}.
The smoothness of $h^{\mu}_{\text{min}}$ is guaranteed because it is a composite function of smooth functions,  $\exp(\cdot)$, summation, and $\ln(\cdot)$, whose domains are properly defined.
The gradient in \eqref{eq:aprxGradMin} is obtained by differentiating \eqref{eq:aprxMin} with respect to $\mathbf{x}$.
\qed

Hereafter, $\mu$ is treated as the parameter associated with the approximation of the minimum value function.

\eqref{eq:aprxMin} and \eqref{eq:aprxGradMin} are the sign-reversed and $\mu$-scaled versions of the LogSumExp function and the softmax function, respectively.
The scaling by $\mu$ plays an important role in the proposed method that enables adjustment between the $\ell_{1}$ norm and the minimum absolute value function.

\subsection{Proposed Penalty: ULPENS}
\label{ssec:ULPENS}

We define the proposed penalty, ULPENS, by combining the smoothly approximated absolute value function in \eqref{eq:abs} and minimum value function in \eqref{eq:aprxMin}, as follows.
\begin{definition}
Let $\mathbf{x}\in\mathbb{R}^{N}$, $\lambda>0$ and $\mu>0$.
ULPENS $\psi:\mathbb{R}^{N}\to\mathbb{R}$ is defined as
\begin{equation}
\label{eq:psi}
    \psi(\mathbf{x}) = -N\mu\ln\left(\sum_{n=1}^{N}\left(\textup{e}^{\frac{x_n}{\lambda}} + \textup{e}^{-\frac{x_n}{\lambda}}\right)^{-\frac{\lambda}{\mu}}\right).
\end{equation}
\end{definition}
The parameters $\lambda>0$ and $\mu>0$ correspond to those in \eqref{eq:abs} and \eqref{eq:aprxMin}, respectively.
ULPENS is an approximation of the $N$-scaled version of \eqref{eq:propId} and is derived by combining \eqref{eq:abs} and \eqref{eq:aprxMin} as follows:
\begin{align}
\label{eq:nmin}
     N\,h^{\mu}_{\min}&\left(\left[h^{\lambda}_{\textup{abs}}(x_1),\dots,h^{\lambda}_{\textup{abs}}(x_N)\right]^{\mathsf{T}}\right) \nonumber \\
     &=-N\mu\ln\left(\sum^{N}_{n=1}  \exp\left(-\frac{1}{\mu}h^{\lambda}_{\textup{abs}}(x_n)\right) \right) \nonumber \\
     &=-N\mu\ln\left(\sum^{N}_{n=1}  \exp\left(-\frac{\lambda}{\mu}\ln\left(\textup{e}^{\frac{x_n}{\lambda}} + \textup{e}^{-\frac{x_n}{\lambda}}\right)\right) \right) \nonumber \\
    &= -N\mu\ln\left(\sum_{n=1}^{N}\left(\textup{e}^{\frac{x_n}{\lambda}} + \textup{e}^{-\frac{x_n}{\lambda}}\right)^{-\frac{\lambda}{\mu}}\right).
\end{align}
The reason of $N$-scaling is given later in Theorem~\ref{thm:L1}.
The gradient of ULPENS is obtained by differentiating \eqref{eq:psi} with respect to $\mathbf{x}$, as follows.
\begin{proposition}
The gradient $\nabla \psi:\mathbb{R}^{N}\to\mathbb{R}^{N}$ is given by
\begin{equation}
\label{eq:gradPsi}
    [\nabla\psi(\mathbf{x})]_{n} = N\varphi_{n}^{(\mathbf{x})}\tanh\left(\frac{x_n}{\lambda}\right),
\end{equation}
where
\begin{equation}
    \varphi_{n}^{(\mathbf{x})} = \left(\sum_{k=1}^{N}\left(\frac{\textup{e}^{\frac{x_n}{\lambda}} + \textup{e}^{-\frac{x_n}{\lambda}}}{\textup{e}^{\frac{x_k}{\lambda}} + \textup{e}^{-\frac{x_k}{\lambda}}}\right)^{\frac{\lambda}{\mu}}\right)^{-1}.
    \label{eq:phi}
\end{equation}
\end{proposition}

\subsection{Approximation Property of ULPENS}

In this subsection, we show that, in a certain sense, ULPENS can be regarded as a function that continuously interpolates between the minimum absolute value function and the $\ell_{1}$ norm.

The following theorem states that ULPENS is an approximation of the $N$-scaled minimum absolute value function.
\begin{theorem}
\label{thm:convergence}
Let $\lambda>0$ and $\mu>0$.
If $\lambda\to+0$ and $\mu\to+0$, then $\psi(\mathbf{x})\to N\min(|\mathbf{x}|)$.
\end{theorem}

\indent\textit{Proof:} The 3rd line of \eqref{eq:nmin} can be rearranged as
\begin{align}
\label{eq:minProof}
    &-N\mu\ln\left(\sum^{N}_{n=1} \exp\left(-\frac{\lambda}{\mu}\ln\left(\textup{e}^{\frac{x_n}{\lambda}} + \textup{e}^{-\frac{x_n}{\lambda}}\right)\right) \right) \nonumber\\
    &=-N\mu\ln\left(\sum^{N}_{n=1}  \exp\left(-\frac{\lambda}{\mu}\ln\left(\textup{e}^{\frac{|x_n|}{\lambda}}\left(1 + \textup{e}^{-2\frac{|x_n|}{\lambda}}\right)\right)\right) \right) \nonumber\\
    &=-N\mu\ln\left(\sum^{N}_{n=1}  \exp\left(-\frac{1}{\mu}\left(|x_n|+\lambda\ln\left(1 + \textup{e}^{-2\frac{|x_n|}{\lambda}}\right)\right)\right) \right)\nonumber\\
    &=-N\mu\ln\left(\textup{e}^{-\frac{x_{\min}(\lambda)}{\mu}}\sum^{N}_{n=1}\textup{e}^{-\frac{z_n(\lambda) - x_{\min}(\lambda)}{\mu}}\right)\nonumber\\
    &= Nx_{\min}(\lambda) - N\mu\ln\left(\sum^{N}_{n=1}\textup{e}^{-\frac{z_n(\lambda) - x_{\min}(\lambda)}{\mu}}\right),
\end{align}
where $z_n(\lambda) = |x_n| + \lambda\ln\left(1 + \textup{e}^{-2\frac{|x_n|}{\lambda}}\right)$, and $x_{\min}(\lambda) = \min\left(\mathbf{z}(\lambda)\right)$.
The index of the minimum value is the same in $\min(\mathbf{z}(\lambda))$ and $\min(|\mathbf{x}|)$ because $|z_n(\lambda)|\ge|z_m(\lambda)|$ when $|x_n|\ge|x_m|$ for all $n, m\in\mathbb{N}$.
Since $\lambda\ln\left(1 + \textup{e}^{-2\frac{|x_n|}{\lambda}}\right)\to+0$ when $\lambda\to+0$, the first term $Nx_{\min}(\lambda)$ in \eqref{eq:minProof} converges to $N\min(|\mathbf{x}|)$ regardless of $\mu$.
The second term in \eqref{eq:minProof} converges to $0$ regardless of the relative convergence speed of $\lambda$ and $\mu$.
To confirm this, we define the continuous functions $\xi$ and $\zeta$ to represent their relationship. 
Specifically, we set $\mu = \xi(\lambda)$, where $\xi(\lambda) \to +0$ as $\lambda \to +0$, and $\lambda = \zeta(\mu)$, where $\zeta(\mu) \to +0$ as $\mu \to +0$.
With these definitions, if $\lambda\to+0$, the second term converges as
\begin{equation*}
    -N\xi(\lambda)\ln\left(\sum^{N}_{n=1}\textup{e}^{-\frac{z_n(\lambda) - x_{\min}(\lambda)}{\xi(\lambda)}}\right) \to 0,
\end{equation*}
and if $\mu\to+0$,
\begin{equation*}
    -N\mu\ln\left(\sum^{N}_{n=1}\textup{e}^{-\frac{z_n(\zeta(\mu)) - x_{\min}(\zeta(\mu))}{\mu}}\right) \to 0,
\end{equation*}
which can be proved in the same manner as Fact~\ref{fact:max}.
\qed

On the $\ell_{1}$-norm side of the approximation, we focus on the gradient behavior of ULPENS.
The following theorem shows that the gradient of ULPENS in \eqref{eq:gradPsi} converges to that of the $\ell_{1}$ norm as $\mu\to+\infty$.

\begin{theorem}
\label{thm:L1}
If $\lambda \to +0$ and $\mu\to+\infty$ then $[\nabla\psi(\mathbf{x})]_n\to\sign(x_n)$.
\end{theorem}
\indent\textit{Proof:} 
We focus on the convergence of $\boldsymbol{\varphi}^{(\mathbf{x})}$ in \eqref{eq:gradPsi}.
$\boldsymbol{\varphi}^{(\mathbf{x})}$ can be rearranged as

\vspace{-9pt}
\footnotesize
\begin{align}
     &\left(\textup{e}^{\frac{x_n}{\lambda}} + \textup{e}^{-\frac{x_n}{\lambda}}\right)^{-\frac{\lambda}{\mu}}\left(\sum_{k=1}^{N}\left(\textup{e}^{\frac{x_k}{\lambda}} + \textup{e}^{-\frac{x_k}{\lambda}}\right)^{-\frac{\lambda}{\mu}}\right)^{-1} \nonumber\\
    &=\textup{e}^{\frac{|x_n|}{\mu}}\left(1 + \textup{e}^{-2\frac{|x_n|}{\lambda}}\right)^{-\frac{\lambda}{\mu}}\left(\sum_{k=1}^{N}\textup{e}^{\frac{|x_k|}{\mu}}\left(1 + \textup{e}^{-2\frac{|x_k|}{\lambda}}\right)^{-\frac{\lambda}{\mu}}\right)^{-1}.\nonumber
\end{align}
\normalsize
Since $\lambda>0$ and $|x_n|\ge 0$, $\left(1 + \textup{e}^{-2\frac{|x_n|}{\lambda}}\right)^{\lambda}\in(1,2^\lambda]$.
When $\mu\to+\infty$, $\left(1 + \textup{e}^{-2\frac{|x_n|}{\lambda}}\right)^{-\frac{\lambda}{\mu}}\to1$ regardless of $\lambda$ because it is $\left(1 + \textup{e}^{-2\frac{|x_n|}{\lambda}}\right)^{\lambda}$, which is a finite value for any $\lambda>0$, raised to the power of $-\frac{1}{\mu}$.
Likewise, $\left(1 + \textup{e}^{-2\frac{|x_k|}{\lambda}}\right)^{-\frac{\lambda}{\mu}}\to1$ as $\mu\to+\infty$.
Since $\textup{e}^{-2\frac{|x_n|}{\mu}}$ and $\textup{e}^{-2\frac{|x_k|}{\mu}}$ converge to 1 as $\mu\to+\infty$, $\boldsymbol{\varphi}^{(\mathbf{x})}$ converges to $1/N$ regardless of $\lambda$.
Owing to the $N$-scaling in \eqref{eq:gradPsi}, this is canceled.
Since $\tanh(\frac{x_n}{\lambda})\to\sign(x_n)$ regardless of $\mu$, $[\nabla\psi(\mathbf{x})]_n\to\sign(x_n)$.
\qed

The following proposition ensures that $\boldsymbol{\varphi}^{(\mathbf{x})}$ in \eqref{eq:gradPsi} continuously changes with respect to $\mu$. 
\begin{proposition}
\label{prop:conti}
When $\lambda>0$ and $\mathbf{x}\in\mathbb{R}^{N}$ are fixed, $\boldsymbol{\varphi}^{(\mathbf{x})}$ is pointwise continuous with respect to $\mu>0$.
\end{proposition}
\indent\textit{Proof:} 
$\boldsymbol{\varphi}^{(\mathbf{x})}$ is the composite function of pointwise continuous functions (with respect to $\mu$), $1/(\cdot)$, summation, and $(\cdot)^{\frac{\lambda}{\mu}}$, whose domain are properly defined.\qed

Note that, although Theorem~\ref{thm:L1} holds for $\nabla\psi$, the function $\psi$ in \eqref{eq:psi} itself cannot be defined when $\mu\to\infty$ because $\psi(\mathbf{x})\to-\infty$ for all $\mathbf{x}$.
On the other hand, Theorem~\ref{thm:convergence} states that if $\mu\to+0$ then $\psi(\mathbf{x})\to N\min(|\mathbf{x}|)$, and Proposition~\ref{prop:conti} guarantees that the gradient changes continuously with respect to $\mu$.
In the sense of these results, ULPENS can be regarded as a function that continuously interpolates between the $\ell_{1}$ norm and the $N$-scaled minimum absolute value function.

\begin{figure}
    \centering
    \includegraphics[width=1\linewidth]{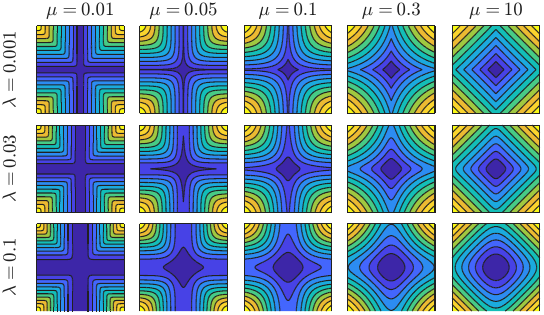}
    \caption{Contour plots of ULPENS in the two-dimensional space.
    $\lambda$ and $\mu$ set for each row and column are shown on the left and top, respectively.
    The parameter $\mu$, which corresponds to the approximation of the minimum function, was set to $0.01, 0.05, 0.1, 0.3, 10$ from left to right.
    The parameter $\lambda$, which corresponds to the approximation of the absolute value function, was set to $0.001, 0.03, 0.1$ from top to bottom.
    These values were chosen for visibility.
    }
    \label{fig:penaFigMany}
\end{figure}

We demonstrate the relationship between parameters and contour plots of ULPENS by Fig.~\ref{fig:penaFigMany}.
In this figure, since the entire contour tends toward $-\infty$ as $\mu$ increases, a bias is added to $\psi(\mathbf{x})$ to counterbalance this effect%
\footnote{The bias is added to each plot so that it has a value of 0 at the origin, i.e., $-\psi(\mathbf{0}) = -2(\mu-\lambda)\ln2$ is added.}.
As $\lambda$ decreases (from bottom to top), the contours near the axes become sharper because the approximation of the absolute value function becomes more accurate.
When $\lambda$ and $\mu$ are small (top left), the contour plot of ULPENS approaches that of the minimum absolute value function, shown in the top left box in Fig.~\ref{fig:pena}.
On the other hand, when $\lambda$ is small and $\mu$ is large (top right), the contour plot of ULPENS approaches that of the $\ell_{1}$ norm, shown in the top right box in Fig.~\ref{fig:pena}.
Note that the functions between each column also exist continuously as ensured by Proposition~\ref{prop:conti}.
These demonstrate that ULPENS enables the continuous transition of the function through the parameter $\mu$ with properly set small $\lambda$.
This behavior will be explained in Sec.~\ref{ssec:weight} through the properties of $\varphi^{(\mathbf{x})}$.

\subsection{ULPENS for Inducing Structured Sparsity}
\label{ssec:structure}

ULPENS can also be used to induce structured sparsity.
We obtain ULPENS for inducing structured sparsity as a composite function of ULPENS in \eqref{eq:psi} and a smoothly approximated $\ell_{2}$ norm, which is defined as follows.
\begin{definition}
\label{def:L2}
Let $\lambda>0$.
The smoothly approximated $\ell_{2}$ norm $h^{\lambda}_{\mathrm{\ell_{2}}}:\mathbb{R}^{N}\to\left[\sqrt{N}(\lambda\ln2),\infty\right)$ is defined as
\begin{equation}
\label{eq:smoothL2}
    h^{\lambda}_{\mathrm{\ell_{2}}}(\mathbf{x}) = \left(\sum^{N}_{n=1}\left(\lambda \ln\left(\textup{e}^{\frac{x_n}{\lambda}} + \textup{e}^{-\frac{x_n}{\lambda}}\right)\right)^{2}\right)^{\frac{1}{2}}.
\end{equation} 
\end{definition}

\begin{proposition}
\label{prop:L2}
$h^{\lambda}_{\mathrm{\ell_{2}}}(\mathbf{x})$ converges to $\|\mathbf{x}\|_{2}$ as $\lambda\to+0$, and is smooth.
Its gradient with respect to $\mathbf{x}$ is given by
\begin{equation}
    [\nabla h^{\lambda}_{\ell_{2}}(\mathbf{x})]_{n} = \tanh\left(\frac{x_n}{\lambda}\right)\cdot\frac{h^{\lambda}_{\textup{abs}}(x_n)}{h^{\lambda}_{\ell_{2}}(\mathbf{x})}.
\end{equation}
\end{proposition}
\indent\textit{Proof:} 
The convergence result follows from \eqref{eq:max}.
The smoothness of $h^{\lambda}_{\mathrm{\ell_{2}}}$ is guaranteed because it is a composite function of smooth functions, $\exp(\cdot)$, summation, $(\cdot)^2$, $\ln(\cdot)$, and $(\cdot)^{\frac{1}{2}}$, whose domains are properly defined.
The gradient is obtained by differentiating \eqref{eq:smoothL2} with respect to $\mathbf{x}$.
\qed

We define ULPENS for inducing structured sparsity by replacing the $\ell_{1}$ norm in the $\ell_{2,1}$ norm in \eqref{eq:L21} with ULPENS in \eqref{eq:psi} and replacing the $\ell_{2}$ norm with $h_{\mathrm{\ell_{2}}}$, as follows.
\begin{definition}
   Let $\mathcal{G}$ be the index sets introduced in Sec.~\ref{ssec:SIP}.
ULPENS for inducing structured sparsity $\mathbf{\psi}^{\mathcal{G}}: \mathbb{R}^{N}\to\mathbb{R}$ is defined as
\begin{equation}
\psi^{\mathcal{G}}(\mathbf{x}) = \psi([h^{\lambda}_{\ell_{2}}(\mathbf{x}_{\mathcal{G}_{1}}),\dots,h^{\lambda}_{\ell_{2}}(\mathbf{x}_{\mathcal{G}_{M}})]^{\mathsf{T}}). 
\end{equation}
\end{definition}

The gradient of $\psi^{\mathcal{G}}$ is obtained by the chain rule as follows.
\begin{proposition}
Let $m(n)\in\mathbb{N}$ be the index of group such that $n\in\mathcal{G}_{m(n)}$.
The $n$th element of the gradient $\nabla\psi^{\mathcal{G}}$ corresponding to $\mathcal{G}_{m(n)}$ is given by
\begin{equation}
\label{eq:gradPsiG}
    [\nabla\psi^{\mathcal{G}}(\mathbf{x})]_{n} = [\nabla \psi(\mathbf{u})]_{m(n)}\cdot [\nabla h^{\lambda}_{\ell_{2}}(\mathbf{x}_{\mathcal{G}_{m(n)}})]_{n}, 
\end{equation}
where $\mathbf{u} = [h^{\lambda}_{\ell_{2}}(\mathbf{x}_{\mathcal{G}_{1}}),\dots,h^{\lambda}_{\ell_{2}}(\mathbf{x}_{\mathcal{G}_{M}})]^{\mathsf{T}} \in \mathbb{R}^{M}$. 
\end{proposition}

\section{Properties and Implementation of ULPENS}
\label{sec:property}

For practical usability of ULPENS, we clarify its properties in terms of ordered weighting and smoothness.
Regarding the ordered weighting, we discuss the characteristics of  $\boldsymbol{\varphi}^{(\mathbf{x})}$ in \eqref{eq:phi} in Sec.~\ref{ssec:weight}.
As for smoothness, we derive an upper bound on the Lipschitz constant of $\nabla\psi$ in Sec.~\ref{ssec:Lip}.
We additionally introduce several implementation tricks in Sec.~\ref{ssec:trick} to stabilize the computation of $\nabla\psi$ and in Sec.~\ref{ssec:tuning} to reduce dependency on the input.

\subsection{Characteristics of $\boldsymbol{\varphi}^{(\mathbf{x})}$ in \eqref{eq:phi}}
\label{ssec:weight}

As in \eqref{eq:gradPsi}, $\boldsymbol{\varphi}^{(\mathbf{x})}$ in \eqref{eq:phi} governs the behavior of the gradient of ULPENS.
The following theorem provides the properties of $\boldsymbol{\varphi}^{(\mathbf{x})}$, which helps explain the property of ULPENS.
\begin{theorem}
\label{thm:phi1}
$\boldsymbol{\varphi}^{(\mathbf{x})}$ in \eqref{eq:phi} has the following properties: 
(a) $\varphi_{n}^{(\mathbf{x})}\in (0,1]$.
(b) $\sum^{N}_{n=1}\varphi^{(\mathbf{x})}_{n} = 1$.
(c) If $|x_n| \le |x_m|$ then $\varphi^{(\mathbf{x})}_{n} \ge \varphi^{(\mathbf{x})}_{m}$.
\end{theorem}

\indent\textit{Proof:}
(b) is shown as

\vspace{-9pt}
\footnotesize
\begin{align*}
    \sum_{n=1}^{N}\varphi_{n}^{(\mathbf{x})} &= \sum_{n=1}^{N}\left(\sum_{k=1}^{N}\left(\frac{\textup{e}^{\frac{x_n}{\lambda}} + \textup{e}^{-\frac{x_n}{\lambda}}}{\textup{e}^{\frac{x_k}{\lambda}} + \textup{e}^{-\frac{x_k}{\lambda}}}\right)^{\frac{\lambda}{\mu}}\right)^{-1} \\
      &= \sum_{n=1}^{N}\left(\left(\textup{e}^{\frac{x_n}{\lambda}} + \textup{e}^{-\frac{x_n}{\lambda}}\right)^{\frac{\lambda}{\mu}}\sum_{k=1}^{N}\left(\textup{e}^{\frac{x_k}{\lambda}} + \textup{e}^{-\frac{x_k}{\lambda}}\right)^{-\frac{\lambda}{\mu}}\right)^{-1}\\
      &= \sum_{n=1}^{N}\left(\textup{e}^{\frac{x_n}{\lambda}} + \textup{e}^{-\frac{x_n}{\lambda}}\right)^{-\frac{\lambda}{\mu}}\left(\sum_{k=1}^{N}\left(\textup{e}^{\frac{x_k}{\lambda}} + \textup{e}^{-\frac{x_k}{\lambda}}\right)^{-\frac{\lambda}{\mu}}\right)^{-1}=1.
\end{align*}
\normalsize
Since $\varphi^{(\mathbf{x})}_n > 0$ and (b) holds, $\varphi^{(\mathbf{x})}_n\in(0,1)$ when $N>1$.
When $N=1$, $\varphi^{(\mathbf{x})}_N = \varphi^{(\mathbf{x})}_1 = 1$.
These yield (a).
(c) can be proved as follows.
$\varphi^{(\mathbf{x})}_n$ can be rearranged so that the denominator is $\sum_{k=1}^{N}\left(\textup{e}^{\frac{x_k}{\lambda}} + \textup{e}^{-\frac{x_k}{\lambda}}\right)^{-\frac{\lambda}{\mu}}$ and the numerator is $\left(\textup{e}^{\frac{x_n}{\lambda}} + \textup{e}^{-\frac{x_n}{\lambda}}\right)^{-\frac{\lambda}{\mu}}$.
The denominators are equal in $\varphi^{(\mathbf{x})}_n$ and $\varphi^{(\mathbf{x})}_m$, and the function of the numerators, $\left(\textup{e}^{\frac{x}{\lambda}} + \textup{e}^{-\frac{x}{\lambda}}\right)^{-\frac{\lambda}{\mu}}$, takes larger value as $|x|$ decreases because it is an even function and monotonically decreases for $x\ge0$.
Hence, (c) holds.
\qed

Theorem~\ref{thm:phi1}(c) means that ULPENS has ordered weighting property.
For inputs with larger magnitudes, $\varphi^{(\mathbf{x})}_n$ takes smaller values, resulting in smaller gradients and minor updates during gradient-based optimization.
Conversely, for inputs with smaller magnitudes, $\varphi^{(\mathbf{x})}_n$ takes larger values, leading to larger gradients and significant updates toward zero.
This behavior demonstrates that ULPENS selectively penalizes smaller values.
From this property, we call $\varphi^{(\mathbf{x})}$ as an \textit{adaptive weight}.

\begin{figure}
    \centering
    \includegraphics[width=1\linewidth]{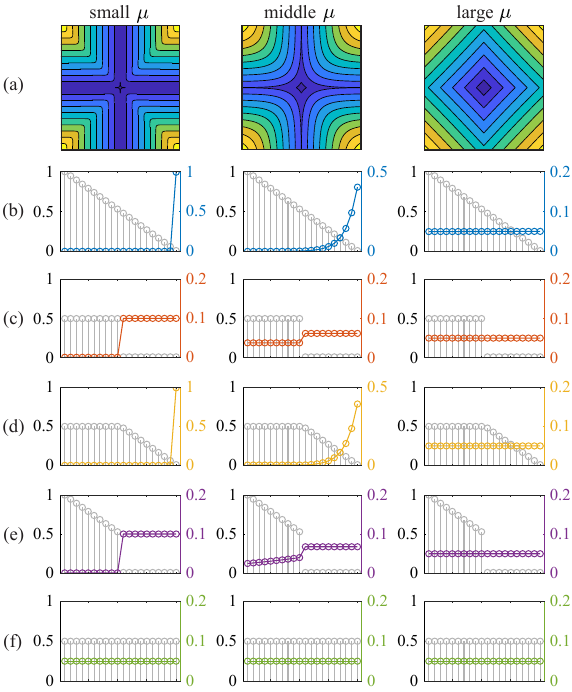}
    \caption{
    Examples of the adaptive weights of ULPENS which change depending on the parameter $\mu$ and input variables.
    The right box has a larger $\mu$.
    (a) shows the contour plot of ULPENS in the two-dimensional space.
    (b--f) show the adaptive weights for some typical 20-dimensional inputs.
    The adaptive weights are colored to have the same color for the same input variable, while the input variables are colored in gray.
    The left and right vertical axes correspond to the inputs and the weight, respectively.
    }
    \label{fig:weight}
\end{figure}

We show an example of adaptive weights $\boldsymbol{\varphi}^{(\mathbf{x})}$ in Fig.~\ref{fig:weight} to illustrate that it has ordered weighting property.
The left column corresponds to a small $\mu$, while the right column corresponds to a large $\mu$.
Figs.~\ref{fig:weight}(b)--(f) show the adaptive weights $\boldsymbol{\varphi}^{(\mathbf{x})}$ (colored) corresponding to 20-dimensional inputs (gray).
A larger weight means the input variable is more penalized.

To check the ordered weighting property, we focus on each columns in Figs.~\ref{fig:weight}(b)--(f).
In (b)--(e) in the middle column (middle $\mu$), the weights corresponding to inputs with smaller absolute values are larger.
These show the ordered weighting property, which enables selective penalization of smaller values.
When $\mu$ is small (left column), in (b) and (d), the adaptive weight corresponding to the smallest input is nearly 1, while the weights for all other inputs are nearly 0.
(c) and (e) are similar in that they assign equally large weights to small inputs and nearly zero weights to large inputs%
\footnote{
The maximum value of the weight depends on the inputs.
If multiple inputs share the same minimum value, the maximum value of the weight is $1/J$, where $J$ is the multiplicity of the minimum value.
In Figs.~\ref{fig:weight}(c) and (e), ten inputs have a minimum value of 0, so the maximum value of the weight is 1/10.
}.
As $\mu$ increases (from left to right), this trend decreases.
When $\mu$ is large (right column), the ordered weighting property disappers and all elements are weighted equally.
From these, $\mu$ determines the extent to which the order of elements is considered.

\subsection{Lipschitz Continuity of $\nabla \psi$ and Its Bound}
\label{ssec:Lip}

For practical use of ULPENS in optimization algorithms, we demonstrate the $L_{\nabla f(\mathbf{x})}$-smoothness of ULPENS.
A function $f$ is $L_{\nabla f(\mathbf{x})}$-smooth in a set $\mathcal{B}$ if the following inequality holds for $L_{\nabla f(\mathbf{x})}\ge0$:
\begin{equation}
\label{eq:defLip}
 \|\mathbf{H}^{(\mathbf{x})}\|_{\text{sp}} \le L_{\nabla f(\mathbf{x})},\quad \forall \mathbf{x}\in\mathcal{B},
\end{equation}
where $\mathbf{H}^{(\mathbf{x})}$ is the Hessian of $f$ at $\mathbf{x}$, and $\|\cdot\|_{\text{sp}}$ is the spectral norm.

We will give a bound for $L_{\nabla\psi(\mathbf{x})}$ in Theorem~\ref{thm:bound}.
To do so, we first calculate the Hessian of ULPENS $\psi$ and prove that  its spectral norm is bounded.
The Hessian is obtained by calculating the second partial derivatives of \eqref{eq:psi}, as follows.

\begin{proposition}
The Hessian of ULPENS at $\mathbf{x}$, denoted as $\mathbf{H}^{(\mathbf{x})} \in \mathbb{R}^{N\times N}$, is given as follows.
The diagonal components are given by
\begin{align}
\label{eq:Jnn}
    H^{(\mathbf{x})}_{n,n}
    &=\frac{N}{\lambda}\varphi_{n}^{(\mathbf{x})} - N\left(\frac{1}{\mu} + \frac{1}{\lambda}\right)\left(\tanh\left(\frac{x_n}{\lambda}\right)\right)^{2}\varphi_{n}^{(\mathbf{x})}\nonumber\\
    &\quad+ \frac{N}{\mu}\left(\tanh\left(\frac{x_n}{\lambda}\right)\right)^{2}\left(\varphi_{n}^{(\mathbf{x})}\right)^{2},
\end{align}
while the non-diagonal components ($n \neq m$) are given by
\begin{equation}
\label{eq:Jmn}
    H^{(\mathbf{x})}_{m,n} = \frac{N}{\mu}\,\tanh\left(\frac{x_n}{\lambda}\right)\,\varphi_{n}^{(\mathbf{x})}\tanh\left(\frac{x_m}{\lambda}\right)\,\varphi_{m}^{(\mathbf{x})},  
\end{equation}
where $\boldsymbol{\varphi}^{(\mathbf{x})}$ is defined in \eqref{eq:phi}.
\end{proposition}

\begin{theorem}
\label{thm:Lip}
$\psi$ is $L_{\nabla\psi(\mathbf{x})}$-smooth.
\end{theorem}
\indent\textit{Proof:} We will show that the eigenvalues of $\mathbf{H}^{(\mathbf{x})}$ is bounded by using the Gershgorin circle theorem \cite{varga}.
To do so, we will separately show the boundedness of diagonal components of $\mathbf{H}^{(\mathbf{x})}$ and radius of the Gershgorin circle calculated by non-diagonal components of $\mathbf{H}^{(\mathbf{x})}$.

When $N=1$, $H^{(\mathbf{x})}_{n,n} = H^{(\mathbf{x})}_{1,1} = \frac{1}{\lambda}(1-\left(\tanh\left(\frac{x_1}{\lambda}\right)\right)^{2})$.
Since $\tanh\left(\frac{x_n}{\lambda}\right)\in(-1,1)$, it is bounded; thus $\psi$ is $L_{\nabla \psi(\mathbf{x})}$-smooth.

We consider the case $N>1$.
The boundedness of the diagonal components $H^{(\mathbf{x})}_{n,n}$ is immediately shown by $\varphi^{(\mathbf{x})}_{n}\in(0,1)$ (in the proof of Theorem~\ref{thm:phi1}(a)) and $\tanh\left(\frac{x_n}{\lambda}\right)\in(-1,1)$.
The boundedness of the radius of the $m$th Gershgorin circle, given by the sum of the absolute values of the non-diagonal elements in the $m$th row, can be shown as follows:
\begin{align}
    \sum_{m}\left|H^{(\mathbf{x})}_{m,n}\right| &=\sum_{m}\left|\frac{N}{\mu}\tanh\left(\frac{x_n}{\lambda}\right)\varphi^{(\mathbf{x})}_{n}\tanh\left(\frac{x_m}{\lambda}\right)\varphi^{(\mathbf{x})}_{m}\right|\nonumber\\
    &= \frac{N}{\mu}\left|\tanh\left(\frac{x_n}{\lambda}\right)\right|\varphi^{(\mathbf{x})}_{n}\sum_{m}\left|\tanh\left(\frac{x_m}{\lambda}\right)\varphi^{(\mathbf{x})}_{m}\right| \nonumber\\
    &< \frac{N}{\mu}\left|\tanh\left(\frac{x_n}{\lambda}\right)\right|\varphi^{(\mathbf{x})}_{n}\sum_{m}\left|\varphi^{(\mathbf{x})}_{m}\right| \nonumber\\
    &= \frac{N}{\mu}\left|\tanh\left(\frac{x_n}{\lambda}\right)\right|\varphi^{(\mathbf{x})}_{n}<\frac{N}{\mu}.
    \label{eq:R}
\end{align}
The 2nd to 3rd line is based on $\tanh\left(\frac{x_n}{\lambda}\right)\in(-1,1)$, and the 3rd to 4th line is based on Theorem~\ref{thm:phi1}(b).

From these, the eigenvalues of $\mathbf{H}^{(\mathbf{x})}$ are bounded because they lie inside the Gershgorin circles, which have centers in the bounded domain and the bounded radius.
Hence $\psi$ is $L_{\nabla \psi}$-smooth.\qed

The above proof provides the following bound of $L_{\nabla\psi(\mathbf{x})}$.
\begin{theorem}
\label{thm:bound}
$L_{\nabla\psi(\mathbf{x})}$ is upper bounded as
\begin{equation}
\label{eq:bound}
    L_{\nabla\psi(\mathbf{x})} < \max(|H^{(\mathbf{x})}_{1,1}|+ \vartheta^{(\mathbf{x})}_{1},\dots,|H^{(\mathbf{x})}_{N,N}|+ \vartheta^{(\mathbf{x})}_{\!N}),
\end{equation}
where $\vartheta^{(\mathbf{x})}_n = \frac{N}{\mu}\left|\tanh\left(\frac{x_n}{\lambda}\right)\right|\varphi^{(\mathbf{x})}_{n}$.    
\end{theorem}
\indent\textit{Proof:}
$\vartheta^{(\mathbf{x})}_1,\dots,\vartheta^{(\mathbf{x})}_N$ are obtained in the final line of \eqref{eq:R} as the bounds of each radius of the Gershgorin circles.
By adding these to the corresponding absolute values of centers and taking the maximum, the upper bound is obtained.
\qed

This bound in \eqref{eq:bound} requires only the computation of the diagonal elements of $\mathbf{H}^{(\mathbf{x})}$, making it computationally less expensive than computing $\mathbf{H}^{(\mathbf{x})}$ and its eigenvalues.
Note that $\frac{N}{\mu}$ can be used instead of $\boldsymbol{\vartheta}^{(\mathbf{x})}$ as shown in \eqref{eq:R}.

\subsection{Implementation Technique for $\nabla \psi$}
\label{ssec:trick}

\begin{algorithm}[t]
{\small
\caption{Calculation of gradient of ULPENS $\nabla\psi(\mathbf{x})$}
\label{alg:diffPena}
\begin{algorithmic}[1]
\renewcommand{\algorithmicrequire}{\textbf{Input:}}
\renewcommand{\algorithmicensure}{\textbf{Output:}}
\REQUIRE $\mathbf{x}\in\mathbb{R}^N,\lambda>0,\mu>0$
\ENSURE $\mathbf{g}$
\STATE $x_{\text{max}} = \max(|\mathbf{x}|)$
\IF{$\mu \ge \lambda$}
\STATE $\mathbf{p} = \left(\text{exp}\left(\frac{|\mathbf{x}|-x_{\text{max}}}{\lambda}\right) + \text{exp}\left(\frac{-|\mathbf{x}|-x_{\text{max}}}{\lambda}\right) \right)^{-\frac{\lambda}{\mu}}$
\ELSE
\STATE $\mathbf{q} = -\frac{\lambda}{\mu}\ln\left( \exp\left(\frac{|\mathbf{x}|-x_{\max}}{\lambda}\right) + \exp\left(\frac{-|\mathbf{x}|-x_{\max}}{\lambda}\right) \right)$
\STATE $q_{\max} = \max(\mathbf{q})$
\STATE $\mathbf{p} = \exp\left(\mathbf{q} - q_{\max}\right)$
\ENDIF
\STATE $\mathbf{y} = \mathbf{p}/\sum^{N}_{n=1}p_n$
\STATE $\mathbf{z} = \left(1 - \text{exp}\left(-\frac{2|\mathbf{x}|}{\lambda}\right)\right)\oslash\left(1 + \text{exp}\left(-\frac{2|\mathbf{x}|}{\lambda}\right)\right)\odot\sign(\mathbf{x})$
\STATE $\mathbf{g} = N\,\mathbf{y}\odot\mathbf{z}$
\end{algorithmic}
}
\end{algorithm}

In this subsection, we introduce several implementation techniques for ULPENS.
The gradient of ULPENS in \eqref{eq:gradPsi}, as well as the hyperbolic tangent function, involves the exponential function, which is prone to overflow.
To address this issue, we implement specific techniques to minimize the risk of overflow.
In addition, the applicable value of $\lambda$ is discussed based on the working precision.

We firstly present the implementation of the gradient of ULPENS in Algorithm~1.
In Algorithm 1, the output $\mathbf{g}$ corresponds to $\nabla \psi$ in \eqref{eq:gradPsi}, $\mathbf{y}$ corresponds to the adaptive weight $\boldsymbol{\varphi}^{(\mathbf{x})}$ in \eqref{eq:phi}, and $\mathbf{z}$ represents the hyperbolic tangent in \eqref{eq:tanh}.  
The operator $\odot$ denotes element-wise multiplication, and $\oslash$ denotes element-wise division.

There are three tricks in Algorithm 1.
The first trick prevents overflow by dividing the numerator and denominator in \eqref{eq:phi} by $\textup{e}^{x_{\max}}$, making the value in $\exp$ less than or equal to 0.
This trick yields the lines 3 and 5 in Algorithm~1. 
The second trick involves calculating the adaptive weight differently depending on whether $\mu\ge\lambda$ (line 3) or $\mu<\lambda$ (lines 5--7).
This is because $(\cdot)^{-\frac{\lambda}{\mu}}$, whose input is smaller than 1 due to the first trick, behaves differently depending on $\frac{\lambda}{\mu}$.
When $\frac{\lambda}{\mu}\le1$, the calculation in line 3 proceeds without issue. 
However, when $\frac{\lambda}{\mu}<1$, this calculation can result in an overflow.
The following procedure in lines 5--7 avoids this issue.
In line 5, the logarithm is taken, and the first trick is applied in lines 6 and 7 to ensure that the input of $\exp$ remains less than or equal to 0. 
Then, in line 7, $\exp$ cancels the logarithm introduced in line 5.
The third trick is to stabilize $\tanh$ in line 10, which prevents the overflow in $\tanh$.
It can be obtained by dividing the numerator and denominator in \eqref{eq:tanh} by $\textup{e}^{|\mathbf{x}|}$, as well as the first trick.

\begin{figure}[t]
    \centering
    \includegraphics[width=1\linewidth]{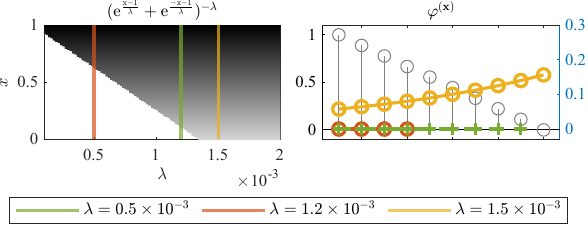}
    \caption{
    The value of $(\textup{e}^{\frac{x-1}{\lambda}} + \textup{e}^{-\frac{-x-1}{\lambda}})^{-\lambda}$ as a function of $x$ and $\lambda$ (left) and the adaptive weight corresponds to the 10-dimensional inputs (gray) computed for each $\lambda$ (right).
    $\lambda$ is set as shown in the legend on the bottom.
    The line colors correspond to each $\lambda$.
    In the left box, larger values appear whiter, and the area where the value is $\infty$ is filled with white.
    In the right box, the left and right vertical axes correspond to the inputs and the weight, respectively.
    Marks are omitted at positions where the weights cannot be computed.
    }
    \label{fig:lambdaTune}
\end{figure}

Next, we discuss how to set $\lambda$.  
In practical use, fine-tuning of $\lambda$ is not necessary; instead, a small value above a certain lower bound can be selected.  
This bound is determined by whether $\mathbf{p}$ and $\mathbf{q}$ in the gradient can be correctly computed%
\footnote{
Although $\mathbf{z}$ also contains $\exp(\cdot)$, it does not cause fatal calculation failure like $\mathbf{p}$ and $\mathbf{q}$ since its values do not diverge.
}, which depends on the working precision.
For explanation, we show an example of $\mathbf{p}$ and the adaptive weight in Fig.~\ref{fig:lambdaTune} ($\mathbf{q}$ is omitted, but the same applies to it).
The inputs are normalized from 0 to 1, $\mu = 1$, and the double-precision arithmetic is employed.
Note that the following explanations are valid for different $\mu$.
As shown, when $\lambda$ is small, the region where the value is $\infty$ (white) appears.
When $\lambda$ is set to such value (red and green), the computation of the adaptive weight fails as shown on the right.
On the other hand, when $\lambda$ is large and $\infty$ does not appear (yellow), the weights can be calculated correctly.
The bound is given by $\frac{\max(|\mathbf{x}|)}{|\ln(c_{\min})|}$, where $c_{\min}$ is the smallest denormal number.
This ensures that the exponential term $\exp\left((|\mathbf{x}| - x_{\text{max}})/\lambda\right)$, which governs the magnitude of the weights in $\mathbf{p}$ and $\mathbf{q}$, remains numerically stable.
For example, in the case of Fig.~\ref{fig:lambdaTune}, the bound is approximately $1.4 \times 10^{-3}$ because $1/|\ln(c_{\min})| \approx 1/|\ln(2.22\times10^{-308})| \approx 1.4 \times 10^{-3}$.
Note that it is more safe to set $c_{\min}$ as the machine epsilon ($2^{-52}$) and choose $\lambda$ such that $\lambda > 1/|\ln(2^{-52})| \approx 2.8 \times 10^{-2}$.

\subsection{Scale-invariant Tuning of $\mu$}
\label{ssec:tuning}

Even when $\mu$ and $\lambda$ are fixed, the adaptive weights in ULPENS vary depending on the maximum absolute value of inputs.
To address this, we propose a tuning method of $\mu$.
The adaptive weight $\boldsymbol{\varphi}^{(\mathbf{x})}$ in \eqref{eq:phi} can be rewritten as
\begin{align}
\label{eq:forScale}
    \varphi^{(\mathbf{x})}_n&=\frac{1}{C^{(\mathbf{x})}}\left(\textup{e}^{\frac{x_n}{\lambda}} + \textup{e}^{-\frac{x_n}{\lambda}}\right)^{-\frac{\lambda}{\mu}}\nonumber\\
    &=\frac{1}{C^{(\mathbf{x})}}\exp\left(\ln\left(\left(\textup{e}^{\frac{x_n}{\lambda}} + \textup{e}^{-\frac{x_n}{\lambda}}\right)^{-\frac{\lambda}{\mu}}\right)\right)\nonumber\\
    &=\frac{1}{C^{(\mathbf{x})}}\exp\left(-\frac{\lambda}{\mu}\ln\left(\textup{e}^{\frac{x_n}{\lambda}} + \textup{e}^{-\frac{x_n}{\lambda}}\right)\right)\nonumber\\
    &=\frac{1}{C^{(\mathbf{x})}}\exp\left(-\frac{h_{\textup{abs}}^{\lambda}(x_{n})}{\mu}\right),
\end{align}
where $C^{(\mathbf{x})}=\sum_{n=1}^N\left(\textup{e}^{\frac{x_n}{\lambda}} + \textup{e}^{-\frac{x_n}{\lambda}}\right)^{-\frac{\lambda}{\mu}}$ is introduced for notational brevity%
\footnote{$C^{(\mathbf{x})}$ does not affect the ordering property of the weight because it just works as a scaling factor.}.
The magnitude of each weight component depends on $\exp(-h^{\lambda}_{\textup{abs}}(x_n)/\mu)$.
Due to the sparsity of the input, the minimum value of $h^{\lambda}_{\textup{abs}}(x_n)$ tends to be close to zero, while its maximum $\phi:=\max(h^{\lambda}_{\textup{abs}}(\mathbf{x}))$ is input-dependent.
As a result, the range $(0,\exp(-\phi/\mu)]$ varies with $\phi$, which in turn alters the shape of the weights.
To avoid this, we normalize $h^{\lambda}_{\textup{abs}}(x_n)$ by its maximum value.
That is, we set $\mu_{\text{tuned}} = \nu\,\phi$, where $\nu > 0$ is a parameter that replaces $\mu$ in controlling the adaptive weight.

\begin{figure}[t]
    \centering
    \includegraphics[width=0.9\linewidth]{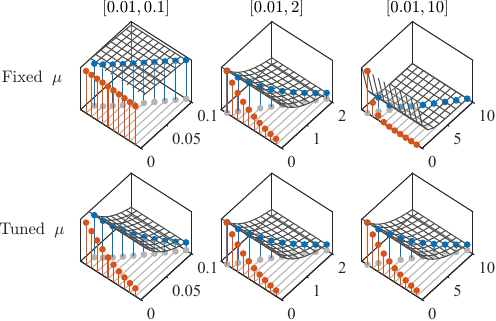}
    \caption{
    Shape of the adaptive weights $\boldsymbol{\varphi}^{(\mathbf{x})}$ for inputs with different maximum values.
    The gray stems represent $h_{\textup{abs}}^{\lambda}(\mathbf{x})$, where $\mathbf{x}$ is a 10-dimensional input with linearly increasing magnitudes.
    The blue stems represent $\exp(-h^{\lambda}_{\textup{abs}}(x_n)/\mu)$, and the black surfaces represent the function $\exp(-\cdot/\mu)$.
    The red stems are projections of the blue stems for improved visibility.
    Note that these stems illustrate the shape of the adaptive weights, since $\boldsymbol{\varphi}^{(\mathbf{x})}$  is obtained by normalizing them with $C^{(\mathbf{x})}$.
    The amplitude range of the input is indicated at the top of each box.
    The top row shows the case where $\mu$ is not tuned, and the bottom row shows the case where $\mu$ is tuned depending on the maximum of input as in Sec.~\ref{ssec:tuning}.
    The parameters are set as $\mu = 1$, $\lambda=0.02$, and $\nu = 0.5$.
    }
    \label{fig:scale}
\end{figure}

An example of the effect of the tuning is shown in Fig.~\ref{fig:scale}. 
As shown in the top boxes, when $\mu$ is fixed, the scale of $\exp(-\cdot/\mu)$ (black surface) varies with the maximum value of the input, resulting in changes in the shape of the weights (red stems) as the maximum input value changes.
In contrast, as shown in the bottom boxes, using $\mu_{\textup{tuned}}$ allows the weight shapes to remain approximately consistent, even when the maximum value of the input varies.

\section{Numerical Experiments}
\label{sec:expts}

\begin{figure*}[t]
    \centering
    \includegraphics[width=0.99\linewidth]{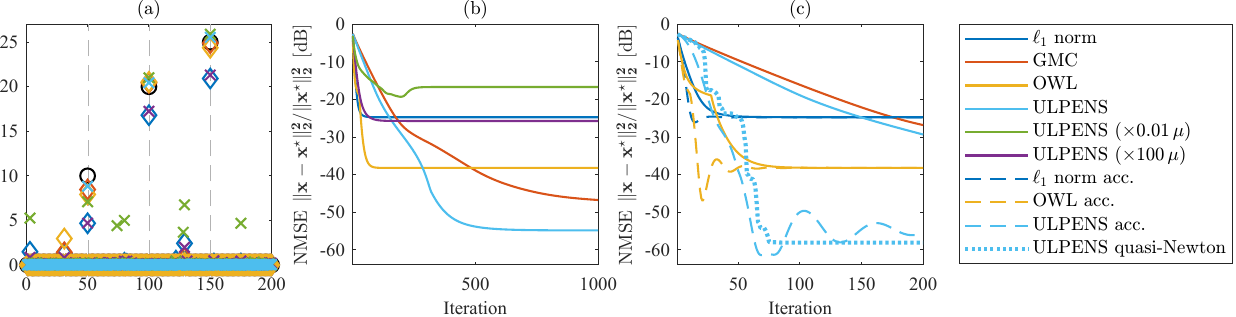}
    \caption{
    Results of the signal recovery in Sec.~\ref{ssec:exptFreq}.
    (a) shows the estimated signals by using each methods.
    (b) shows the NMSE ($\|\mathbf{x}-\mathbf{x}^{\star}\|^{2}_{2}/\|\mathbf{x}^{\star}\|^{2}_{2}$) per iteration without acceleration, where $\mathbf{x}^{\star}$ is the ground truth.
    (c) shows the NMSE per iteration with the Nesterov acceleration and the quasi-Newton method.
    The legends for the figures are shown in the right box, where the dashed lines represent the methods using the Nesterov acceleration, and the dotted lines represent that using the quasi-Newton method.
    The colors of the methods in (a) correspond to those in (b) and (c).  
    For clarity, the ground truth are shown as black circle and its location as gray dashed vertical lines, the conventional methods as diamonds, and the proposed method as crosses.
    }
    \label{fig:freq1}
\end{figure*}

To evaluate the proposed method, we conducted four experiments.
We applied ULPENS to three signal recovery problems and compared with conventional penalties.
We also evaluated the tightness of the upper bound of $L_{\nabla\psi(\mathbf{x})}$ in \eqref{eq:bound}.

\subsection{Sparsity-based Signal Recovery}
\label{ssec:exptFreq}

In this experiment, the performance of sparse signal recovery was compared for different regularization functions and acceleration methods.
We considered the following setting and formulated the corresponding optimization problem.
Let a signal $\mathbf{s} \in \mathbb{R}^{300}$ be observed as $\mathbf{s} = \mathbf{A}\mathbf{x}^{\star} + \mathbf{n}$, where $\mathbf{A}\in \mathbb{R}^{150\times300}$, $\mathbf{x}^{\star}\in\mathbb{R}^{300}$ is the sparse signal to be estimated, and $\mathbf{n}\in\mathbb{R}^{150}$ is additive white Gaussian noise.
The SNR of the observed signal was 5~dB.
The ground truth $\mathbf{x}^{\star}$ was defined as $x^{\star}_{50} = 10$, $x^{\star}_{15} = 20$, and $x^{\star}_{150} = 25$, with all other elements were zero, as indicated by the black circles in Fig.~\ref{fig:freq1}(a).
$\mathbf{A}$ was generated by first creating a matrix with elements drawn from a uniform distribution over $[0,1]$ and then replacing its singular values with values linearly varying from $10^{-4}$ to $1$.
The optimization problem for estimating the signal was formulated as finding $\mathbf{x}\in\mathbb{R}^{300}$ that minimizes \eqref{eq:cost}, where the data fidelity term was $g_{\mathbf{s}}(\mathbf{x}) = \frac{1}{2}\|\mathbf{A}\mathbf{x}-\mathbf{s}\|^{2}_{2}$.

For the regularization term, we used the $\ell_{1}$ norm, GMC \cite{Ivan2017}, OWL \cite{sasaki2024} and ULPENS.
For the conventional method with the $\ell_{1}$ norm and OWL, we used PGM with the step size parameter was fixed at $1$ for all iterations, since the Lipschitz constant of the function $\frac{1}{2}\|\mathbf{A}\mathbf{x} - \mathbf{s}\|_{2}^{2}$ was 1.
For the conventional method with GMC, we used the algorithm from the original paper \cite{Ivan2017}.
For the proposed method, we used GD with the step size parameter was set to $1/(1+\gamma\,L_{\nabla\psi(\mathbf{x})})$ for each iteration.
This setting makes the comparison fair because all step sizes were set to one over the Lipschitz constant.
We also compared the performance of the accelerated methods at each iteration.
For both the conventional penalties and ULPENS, we employed Nesterov acceleration.
Additionally, since only the proposed method is smooth, we applied the quasi-Newton method to it.
For the quasi-Newton method, we used the $\texttt{fmincon}$ function in MATLAB with a termination condition that the norm of the gradient was $10^{-6}$, employing the limited-memory BFGS method.
For evaluation, the normalized mean square error (NMSE) was used.
The initial value of $\mathbf{x}$ was set to the observation $\mathbf{s}$.
The regularization parameter for all methods, the weight parameter for OWL, and ULPENS paramters $\lambda$ and $\mu$ were tuned to achieve the best performance by using the $\texttt{bayesopt}$ function in MATLAB.

\subsubsection{Comparison with the Conventional Methods}
The result is shown in Fig.~\ref{fig:freq1}.
We first compare the methods without acceleration.
As shown in Fig.~\ref{fig:freq1}(b), the proposed method showed the lowest NMSE among all. 
Fig.~\ref{fig:freq1}(a) further demonstrates that ULPENS effectively prevented the underestimation problem compared to the $\ell_{1}$ norm.
In addition, while GMC and OWL failed to suppress some components, ULPENS successfully reduced it.
To show the dependency on $\mu$ of ULPENS, we also used $\mu$ scaled by factors of $0.01$ and $100$ relative to the value obtained via \texttt{bayesopt}.
When $\mu$ was large ($\times 100$), the behavior of ULPENS closely resembled that of the $\ell_{1}$ norm, as indicated by the purple crosses and lines in Figs.~\ref{fig:freq1}(a) and (b).
On the other hand, when $\mu$ was small ($\times 0.01$), some components was not suppressed, as shown by the green crosses in Fig.~\ref{fig:freq1}(a), resulting in a higher NMSE, as shown by the green line in Fig.~\ref{fig:freq1}(b).
This is because reducing $\mu$ causes ULPENS to approximate the minimum absolute value function, which selectively suppresses only a limited number of small values while leaving the others unchanged.

Next, we compare the performance of the accelerated methods.
As shown in Fig.~\ref{fig:freq1}(c), Nesterov acceleration (dashed lines) led to faster convergence for all methods compared to the original methods (solid lines). 
In particular, OWL (yellow dashed line) converged in approximately 100 iterations.
We compared these results with those of the proposed method using the quasi-Newton method.
The proposed method (light blue dotted line) converged in about 70 iterations, achieving the fastest convergence and the best performance.
These results confirm that the smoothness of ULPENS enables efficient optimization using the quasi-Newton method.

\begin{figure*}[t]
    \centering
    \includegraphics[width=0.99\linewidth]{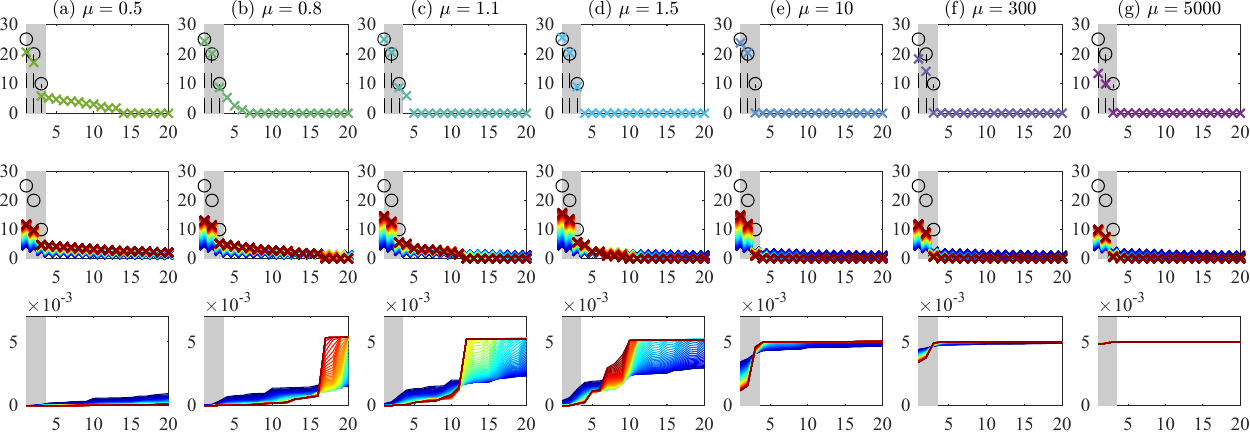}
    \caption{
    Illustration of the first 20 sorted estimated signals and their corresponding adaptive weights evolving with iteration for different $\mu$ settings (Sec.~\ref{ssec:exptFreq}2).
    The top row shows the estimated signals at the final iteration (1000th iteration), the middle row shows the estimated signals from the first to the 100th iteration, and the bottom row shows the corresponding adaptive weights over the same iterations.
    Black circles in the top and middle rows are the ground truth, and the gray areas are their locations.
    The righter the columns, $\mu$ was set larger as shown on the top of each column.
    These $\mu$ values were chosen for visibility.
    The estimated signals and adaptive weights in the middle and bottom rows are colored according to the iteration number.
    Blue colors represent the early iterations, while red colors correspond to iterations closer to 100.
    }
    \label{fig:weightVec}
\end{figure*}

\begin{figure}
    \centering
    \includegraphics[width=1\linewidth]{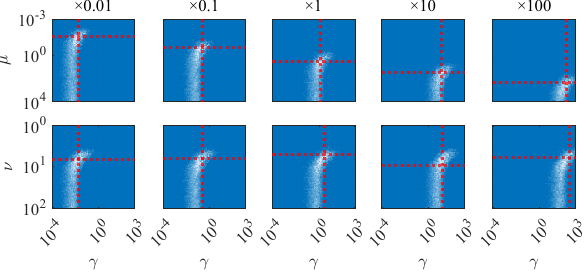}
    \caption{
    Illustration of the NMSE at the final iteration under varying input scales (Sec.~\ref{ssec:exptFreq}3).
    The top row shows the results obtained using $\mu$, while the bottom row presents the results using $\nu$ in the scale-invariant tuning method (Sec.~\ref{ssec:tuning}).
    Each column corresponds to different input scales, which are indicated above the columns.
    Brighter (whiter) pixels indicate lower NMSE values.
    The combination of $\mu$, $\nu$ and $\gamma$ that achieves the best performance under each scale condition is marked by red dotted lines.}
    \label{fig:changeMu}
\end{figure}

\subsubsection{Changes of the Adaptive Weight}
Fig.~\ref{fig:weightVec} illustrates the evolution of the adaptive weights over iterations under different $\mu$.  
This illustration corresponds to those in Fig.~\ref{fig:weight}, but the weights were plotted by lines for better visibility.
The closer a weight was to zero, the smaller its gradient became, meaning that the corresponding component was less suppressed.
When $\mu = 1.5$, the three non-zero elements were correctly estimated, as shown in the center column of Fig.~\ref{fig:weightVec}(d).  
In this case, the weights were gradually adjusted throughout the iterations to retain the three non-zero components while suppressing the remaining zero components.
When $\mu$ was small, as illustrated in Figs.~\ref{fig:weightVec}(a), (b) and (c), the weights tended to retain many components.
This behavior reflects the property of ULPENS approaching the minimum value function, yielding weights similar to those shown in the left column of Fig.~\ref{fig:weight}.
As $\mu$ increased to $1.5$, the weights changed to suppress more components, as shown in the bottom boxes in Figs.~\ref{fig:weightVec}(a), (b) and (c).
As $\mu$ increased beyond 1.5 (Figs.~\ref{fig:weightVec}(e) and (f)), the weights changed to suppress larger components.
Finally, when $\mu$ was even larger, as in Fig.~\ref{fig:weightVec}(g), the weights became uniform, corresponding to the $\ell_{1}$ norm.
With appropriately chosen parameters, the weights can adapt to preserve the components that should be retained, which may lead to improved performance.

\subsubsection{Effect of the Scale-invariant Tuning in Sec.~\ref{ssec:tuning}}

To show the effect of the tuning, we investigated the appropriate settings of $\mu$ and $\nu$ when the scale of the ground truth signal $\mathbf{x}^{\star}$ was changed.
Each of $\mu$, $\nu$, and $\gamma$ was assigned 100 logarithmically spaced values over $[10^{-4}, 10^{-3}]$, $[10^0, 10^2]$, and $[10^{-4}, 10^3]$, respectively.
The result is shown in Fig.~\ref{fig:changeMu}.
Without scale-invariant tuning, the appropriate value of $\mu$ depends on the scale of the signal, as shown by the horizontal red dotted line in the top row of  Fig.~\ref{fig:changeMu}.
On the other hand, with scale-invariant tuning, the appropriate value of $\nu$ was almost constant regardless of the scale of the signal.
This result confirms that the use of scale-invariant tuning can simplify the tuning process.

\subsection{Signal Recovery Based on Group Sparsity}
\label{ssec:exptGroup}

\begin{figure*}
    \centering
    \includegraphics[width=1\linewidth]{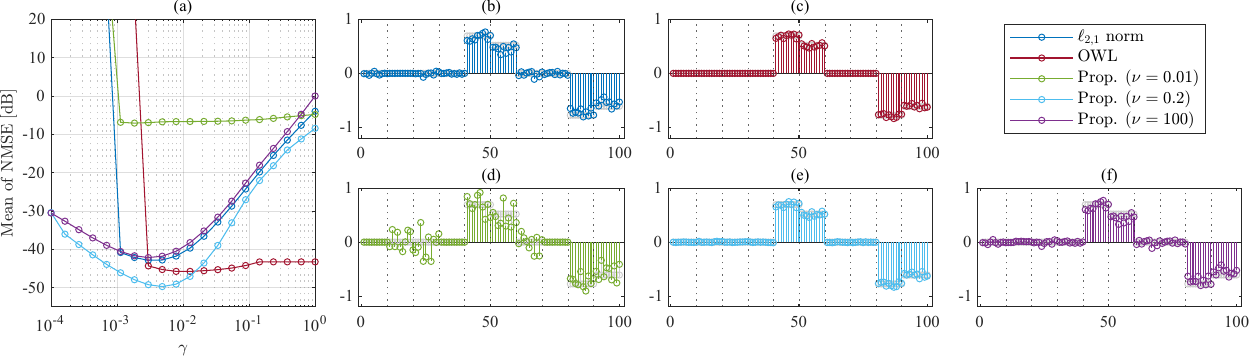}
    \caption{
    Results of signal recovery based on group sparsity in Sec.~\ref{ssec:exptGroup}.
    (a) shows the mean NMSE of 100 trials for each method and each regularization parameter $\gamma$.
    The example of recovered signal are shown in (b)--(f).
    The color of each plot follows the legend on the right top.
    In (b)--(f), the dotted vertical lines indicate the boundaries of the group, and the gray stems indicate the ground truth $\mathbf{x}^{\star}$. 
    }
    \label{fig:structure}
\end{figure*}

In this experiment, to show the effectiveness of ULPENS that induces structured sparsity (Sec.~\ref{ssec:structure}), we conducted a signal recovery experiment based on the group sparsity of block structured signal.
We considered the following setting and formulated the corresponding optimization problem.
We assumed that the measurement was given by $\mathbf{s} = \mathbf{A}\mathbf{x}^{\star} + \mathbf{n}$, where $\mathbf{n} \in \mathbb{R}^{75}$ represents additive white Gaussian noise.
The model matrix $\mathbf{A}\in\mathbb{R}^{75\times 100}$ was generated in the same procedure as in Sec.~\ref{ssec:exptFreq}.
The SNR of the observed signal was 40\,dB.
A structured signal $\mathbf{x}^{\star}\in\mathbb{R}^{100}$ with 10 blocks (10 elements per block) was generated, where 4 blocks were randomly selected to contain nonzero values drawn from $\mathcal{N}(0,1)$, as shown by the gray stems in Fig.~\ref{fig:structure}.
The group of indices $\mathcal{G}_{1},\dots\mathcal{G}_{10}$ that had this structure was defined as $\mathcal{G}_{m} = \{10(m-1)+1,\dots,10m\}\;(1\le m\le 10)$.
The optimization problem for recovering the structured signal was formulated as finding $\mathbf{x} \in \mathbb{R}^{100}$ that minimizes \eqref{eq:cost}, where the data fidelity term was $g_{\mathbf{s}}(\mathbf{x}) = \frac{1}{2}\|\mathbf{A}\mathbf{x} - \mathbf{s}\|_2^2$, and the sparse regularization term was $r(\mathbf{x}) = r^{\mathcal{G}}(\mathbf{x})$, which promotes structured sparsity corresponding to the group of indices $\mathcal{G} = \{\mathcal{G}_{1},\dots,\mathcal{G}_{10}\}$.

\begin{figure*}
    \centering
    \includegraphics[width=1\linewidth]{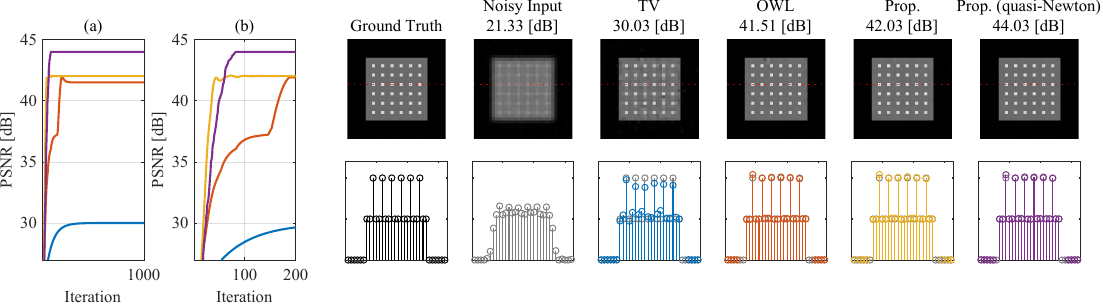}
    \caption{
    PSNR and the corresponding deblurred images in Sec.~\ref{ssec:exptSmoothing}.
    (a) shows PSNR for each iteration, and (b) shows its enlarged part of early 200 iterations.
    Excluding (a) and (b), the top row presents the ground truth, the noisy input, and the deblurred images from left to right.
    The bottom row highlights a portion of the image along the red dotted line.
    In the bottom right four boxes, the ground truth is indicated by the gray stems.
    The line colors in these boxes correspond to (a) and (b): blue indicates TV, red indicates OWL, yellow indicates the proposed method, and purple indicates the proposed method using the quasi-Newton method.
    PSNRs at the final iteration corresponding to each image are shown above the images.
    }
    \label{fig:smooth}
\end{figure*}

We compared the proposed method to the $\ell_{2,1}$ norm and OWL extended to handle the structured sparsity \cite{sasaki2024}.
For each method, the regularization parameter $\gamma$ was varied from $10^{-3}$ to $10$, and the mean of NMSE over 100 trials was compared.
The weights of OWL were tuned to obtain the best performance.
To solve the problem using each penalty function, we used PGM for the conventional methods and GD for the proposed method.
In the PGM, the step size parameter was set to $1$.
For GD, the step size parameter was set to $1/(1+\gamma\,L_{\nabla\psi(\mathbf{x})})$ for each iteration.
In the proposed method, we set $\lambda = 0.1$ and $\nu = 0.01, 0.2, 100$.
The iteration was terminated when either the norm of the difference between successive iterates fell below the threshold of $10^{-7}$ or the maximum iteration count of $10^5$ was reached.

The results are shown in Fig.~\ref{fig:structure}.
The proposed method achieved the best performance with the lowest average NMSE, when $\nu=0.2$ and $\gamma$ was around $3\times10^{-3}$, as shown in the light blue line in Fig.~\ref{fig:structure}(a).
On the other hand, as seen in Fig.~\ref{fig:structure}(d), although the proposed method reduced some blocks to zero when $\nu=0.01$, nonzero values remained in blocks that should be zero, such as $\mathcal{G}_{2}$, $\mathcal{G}_{3}$, and $\mathcal{G}_{7}$, resulting in a high NMSE.

We also examine the performance when $\gamma$ is small, as shown in Fig.~\ref{fig:structure}(a).
For sufficiently large $\nu$, unlike conventional methods, ULPENS does not exhibit extremely high NMSE even for small $\gamma$ ($<10^{-3}$), as indicated by the light blue and purple lines in Fig.~\ref{fig:structure}(a).
This may be due to the smoothness of ULPENS. 
On the other hand, as shown by the green line in Fig.~\ref{fig:structure}(a), when $\nu$ is small, ULPENS exhibited extremely high NMSE.
This may be because ULPENS is close to the minimum absolute value function, which is less smooth.
These results suggest that ULPENS has the advantage of enabling adjustment of the degree of non-convexity.

\subsection{Application to Image Deblurring}
\label{ssec:exptSmoothing}

In this experiment, we applied ULPENS that induces structured sparsity to image deblurring.
We considered the following setting and formulated the corresponding optimization problem.
We blurred the ground truth image shown in Fig.~\ref{fig:smooth}, whose pixel values were 1.0, 0.5 or 0.
The observation $\mathbf{s}\in[0,1]^{N^2}$ was given as $\mathbf{s} = \mathbf{A}\mathbf{x} + \mathbf{n}$ with $N=32$, where $\mathbf{A}\in[0,1]^{N^2\times N^2}$ represents the two-dimensional Gaussian blur matrix with its standard deviation $\sigma=0.75$, $\mathbf{n}$ is additive white Gaussian noise.
The maximum and minimum singular values of $\mathbf{A}$ was 1 and $1.55\times10^{-2}$, respectively.
The optimization problem for recovering a deblurred image was formulated as finding $\mathbf{x}\in[0,1]^{N^2}$ that minimizes \eqref{eq:cost}, where the data fidelity term was $g_{\mathbf{s}}(\mathbf{x}) = \frac{1}{2}\|\mathbf{A}\mathbf{x}-\mathbf{s}\|^{2}_{2}$, and the sparse regularization term that promotes smoothness was $r(\mathbf{x}) = r^{\mathcal{G}}(\mathbf{D}\mathbf{x})$.
Here, $\mathbf{D}\in\{-1,0,1\}^{2N(N-1) \times N^2}$ was difference matrix, and $\mathcal{G} = \{\mathcal{G}_{1},\dots,\mathcal{G}_{N(N-1)}\}$ was a set of groups of indices that combine the vertical and horizontal differences at each pixel.

We compared the proposed method with the $\ell_{2,1}$ norm, which results in the total variation (TV) regularization, and a method \cite{sasaki2024} in which the $\ell_{2,1}$ norm in TV is replaced with OWL.
To solve the problem of using each penalty function, we used ADMM for the conventional methods and GD with Nesterov acceleration and the quasi-Newton method for the proposed method. 
In the gradient method, the step size parameter was set to $1/(1+\gamma\,L_{\nabla\psi(\mathbf{x})})$ for each iteration.
For the quasi-Newton method, we used the $\texttt{fmincon}$ function in MATLAB as described in Sec.~\ref{ssec:exptFreq}.
The regularization parameter for all methods, the weight parameter for OWL, and ULPENS paramters $\lambda$ and $\nu$ were tuned to achieve the best performance by using the $\texttt{bayesopt}$ function in MATLAB.
The denoising results were evaluated by PSNR.

Fig.~\ref{fig:smooth} shows the results of image deblurring.
As shown in Figs.~\ref{fig:smooth}(a) and (b), the proposed method (yellow and purple) achieved a higher PSNR than the conventional methods.
Additionally, as seen in the bottom-right boxes of Fig.~\ref{fig:smooth}, compared to TV regularization, which excessively smoothed the fine structures of the image, OWL and the proposed method effectively smoothed the image while preserving fine details.
Furthermore, the proposed method using the quasi-Newton method converged in only 80 iterations and showed the highest PSNR, as shown in Fig.~\ref{fig:smooth}(b).

\subsection{$L_{\nabla \psi(\mathbf{x})}$ and Tightness of the Bound on $L_{\nabla \psi(\mathbf{x})}$}
\label{ssec:exptLips}

\begin{figure*}
    \centering
    \includegraphics[width=0.98\linewidth]{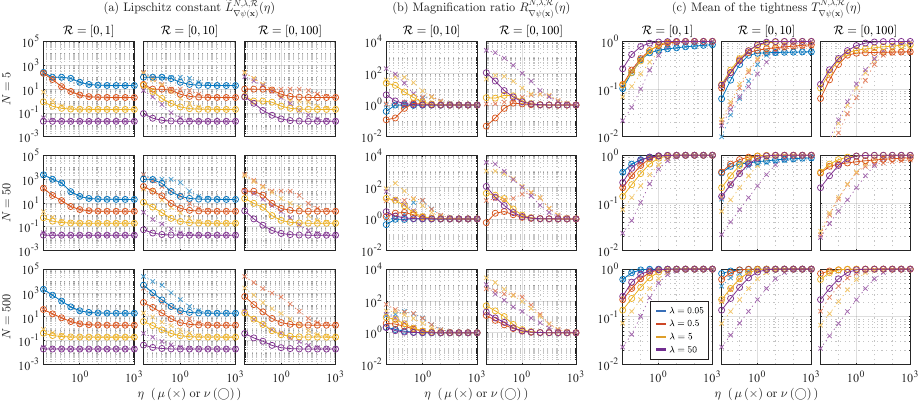}
    \caption{
    The Lipschitz constant of $\nabla\psi$ and the tightness of the bounds calculated using 10\,000 trials of $\mathbf{x}$ (Sec.~\ref{ssec:exptLips}).
    (a) The numerically calculated Lipschitz constant $\widetilde{L}^{N,\lambda,\mathcal{R}}_{\nabla\psi(\mathbf{x})}(\eta)$ and (c) mean of the tightness $T^{N,\lambda,\mathcal{R}}_{\nabla\psi(\mathbf{x})}(\eta)$.
    (b) shows the magnification ratios corresponding to input ranges $\mathcal{R}=[0,10]$ and $\mathcal{R}=[0,100]$ (left: $R^{N,\lambda,[0,10]}_{\nabla\psi(\mathbf{x})}(\eta)$, right: $R^{N,\lambda,[0,100]}_{\nabla\psi(\mathbf{x})}(\eta)$).
    In each box, crosses and dotted lines denote the results for $\mu$, whereas circles and solid lines denote the results for $\nu$.
    The blue, red, yellow and purple lines indicates $\lambda = 0.05, 0.5, 5, 50$, respectively.
    The blue lines were missing in the right column of each diagram because they could not be calculated, as their settings fell below the lower bound of $\lambda$ described in Sec.~\ref{ssec:trick}.
    To illustrate the behavior across different $\lambda$, we selected the four $\lambda$ values.
    }
    \label{fig:Lips}
\end{figure*}

Finally, we computed the Lipschitz constant of $\nabla\psi$ numerically and examined the tightness of the bound given in \eqref{eq:bound}.
We generated $10\,000$ random inputs of varying sizes ($N = 5, 50, 500$) and magnitudes, drawn from uniform distributions over $[0,1]$, $[0,10]$ and $[0,100]$. 
The range of the distributions is denoted by $\mathcal{R}$, i.e., the input range.
We numerically calculated the Lipschitz constant $\widetilde{L}^{N,\lambda,\mathcal{R}}_{\nabla\psi(\mathbf{x})}(\eta)$ as the maximum absolute eigenvalue of \(\mathbf{H}^{(\mathbf{x})}\), evaluated over \(10\,000\) trials of $\mathbf{x}$.
In the notation $\widetilde{L}^{N,\lambda,\mathcal{R}}_{\nabla\psi(\mathbf{x})}(\eta)$, the superscript denotes the parameters: \(\mathcal{R}\) denotes the magnitude range of the input, and \(\eta\) corresponds to \(\mu\) or \(\nu\).
To evaluate the tightness of the bounds, we introduce the following values as \textit{tightness}:
\begin{equation}
    T^{N,\lambda,\mathcal{R}}_{\nabla\psi(\mathbf{x})}(\eta) = \left.\widetilde{L}^{N,\lambda,\mathcal{R}}_{\nabla\psi(\mathbf{x})}(\eta) \middle/ B^{N,\lambda}_{\nabla\psi(\mathbf{x})}(\eta)\right.,
\end{equation}
where $B^{\lambda,N}_{\nabla\psi(\mathbf{x})}(\eta)$ is the bound (the right-hand side of \eqref{eq:bound}). 
This ratio serves as a tightness measure, the closer this value is to 1, meaning that the inequality in \eqref{eq:bound} is tighter.
We used the mean tightness over the $10\,000$ trials of $\mathbf{x}$ for the evaluation.

First, we discuss the Lipschitz constant of $\nabla\psi$. 
Fig.~\ref{fig:Lips}(a) shows $\widetilde{L}^{N,\lambda,\mathcal{R}}_{\nabla\psi(\mathbf{x})}(\eta)$ for each condition.
As $\mu$ and $\nu$ increase, the colored curves in Fig.~\ref{fig:Lips}(a) converge to constant values at the right end for each $\lambda$.
These constants are in fact equal to $1/\lambda$.
This is because, as $\mu \to \infty$, the gradient of ULPENS converges to the element-wise $\tanh(\cdot/\lambda)$ function, whose Lipschitz constant is $1/\lambda$, since $\varphi_{n}^{(\mathbf{x})} = 1/N$ in \eqref{eq:gradPsi}.
As $\mu$ and $\nu$ decreased, $\widetilde{L}^{N,\lambda,\mathcal{R}}_{\nabla\psi(\mathbf{x})}(\eta)$ tended to increase, possibly due to ULPENS approaching the minimum absolute value function.
We also compare the behavior of $\widetilde{L}^{N,\lambda,\mathcal{R}}_{\nabla\psi(\mathbf{x})}(\eta)$ with respect to $\mu$ and $\nu$ for different input ranges. To this end, we define the \textit{magnification ratio} as
\begin{equation}
    R^{N,\lambda,\mathcal{R}}_{\nabla\psi(\mathbf{x})}(\eta) =
    \left.\widetilde{L}^{N,\lambda,\mathcal{R}}_{\nabla\psi(\mathbf{x})}(\eta)
    \middle/
    \widetilde{L}^{N,\lambda,[0,1]}_{\nabla\psi(\mathbf{x})}(\eta)\right.,
\end{equation}
where a value closer to $1$ indicates that the Lipschitz constant changes little from $\widetilde{L}^{N,\lambda,[0,1]}_{\nabla\psi(\mathbf{x})}(\eta)$ when the input range changes.
By comparing the markers and lines of the same color in Fig.~\ref{fig:Lips}(b), for $\nu$ (circle and solid), both $R^{N,\lambda,[0,10]}_{\nabla\psi(\mathbf{x})}(\eta)$ and $R^{N,\lambda,[0,100]}_{\nabla\psi(\mathbf{x})}(\eta)$ were closer to $1$ compared to the case of $\mu$ (cross and dot), especially when $\nu$ was not small.
This indicates that the Lipschitz constant for $\nu$ remained nearly the same across different input ranges.

Finally, the proposed bound is examined.
Fig.~\ref{fig:Lips}(c) shows the mean tightness for each condition.
As indicated by all the boxes, when $\nu$ is relatively large, the bound was tight, whereas for $\nu$ below about 1, the mean of the tightness decreases under all the conditions.
In particular, at $\nu=10^{-2}$, it reaches values in the range of about 0.01--0.1.
In addition, for $\mu$, the tightness of the bound varied with changes in the input range, while for $\nu$, it remained almost unchanged.
These results indicate that the bound becomes less tight when $\nu$ is small, and that incorporating $\nu$ makes the bound more robust to variations in the input range.

\section{Conclusion}
\label{sec:conclusion}

In this paper, we proposed ULPENS, a non-convex, non-separable and smooth penalty function inducing sparsity, that addresses the underestimation problem.
By applying the ultra-discretization formula, we obtained ULPENS in \eqref{eq:psi} as a function that continuously interpolates between the minimum absolute value function and the $\ell_{1}$ norm, achieving selective suppression.
This behavior arises from the adaptive weight in \eqref{eq:phi} that has the ordered weighting property.
We demonstrated the smoothness of ULPENS and the upper bound on the Lipschitz constant of the gradient for practical usability.
We also introduced the implementation technique and effective tuning method.
Our experiments confirmed the effectiveness and applicability of the ULPENS to sparse optimization problems, such as denoising and image deblurring.
In particular, high-performance and fast optimization using the quasi-Newton method was achieved thanks to the selective suppression and smoothness of ULPENS.
The performance, however, depends on appropriate adjustment of the parameter $\nu$.
In future work, we aim to investigate automatic tuning method driven by incoming signals.

\section*{Acknowledgment}
\noindent This work was partly supported by JST CREST, Japan, Grant Number JPMJCR2233.

\bibliographystyle{IEEEtran}
\bibliography{IEEEabrv,Bibliography}

\end{document}